\documentclass[11pt]{amsart}
\usepackage{latexsym,amssymb,amsmath}
 2

\newtheorem{thm}{Theorem}[section]

\newtheorem{prop}[thm]{Proposition}

\newcommand{\thmref}[1]{Theorem~\ref{#1}}

\newcommand{\propref}[1]{Proposition~\ref{#1}}

\begin{document}

\title[Quadratic forms in eight variables]
{On the representations of a positive integer by certain classes of quadratic forms in eight  variables}
\author{B. Ramakrishnan, Brundaban Sahu and Anup Kumar Singh}
\address[B. Ramakrishnan and Anup Kumar Singh]{Harish-Chandra Research Institute, 
       Chhatnag Road, Jhunsi,
     Allahabad -     211 019,
   India.}
\address[Brundaban Sahu]
{School of Mathematical Sciences, National Institute of Science 
Education and Research,
PO: Sainik School,
Bhubaneswar, Odisha - 751 005,
India.}

\email[B. Ramakrishnan]{ramki@hri.res.in}
\email[Brundaban Sahu]{brundaban.sahu@niser.ac.in}
\email[Anup Kumar Singh]{anupsingh@hri.res.in}
%\dedicatory{Dedicated to Professor Krishnaswamy Alladi on the occasion of his 60th birthday}

%\date{July 16, 2016}

\begin{abstract}
In this paper we use the theory of modular forms to find formulas for the number of 
representations of a positive integer by certain class of quadratic forms in eight variables, viz., forms of the form $a_1x_1^2 + a_2 x_2^2 + a_3 x_3^2 + a_4 x_4^2 + b_1(x_5^2+x_5x_6 + x_6^2) + b_2(x_7^2+x_7x_8 + x_8^2)$, where $a_1\le a_2\le a_3\le a_4$, $b_1\le b_2$ and  $a_i$'s $\in \{1,2,3\}$,  $b_i$'s $\in \{1,2,4\}$. We also determine formulas for the number of representations of a positive integer by the quadratic forms $(x_1^2+x_1x_2+x_2^2) + c_1(x_3^2+x_3x_4+x_4^2) + c_2(x_5^2+x_5x_6+x_6^2) + c_3(x_7^2+x_7x_8+x_8^2)$,
where $c_1,c_2,c_3\in \{1,2,4,8\}$, $c_1\le c_2\le c_3$.. 
\end{abstract}

\maketitle

%\bigskip

\section{Introduction}
In this paper we consider the problem of finding the number of representations of the following quadratic forms in eight variables given by 
\begin{equation}\label{quad1}
a_1x_1^2 + a_2 x_2^2 + a_3 x_3^2 + a_4 x_4^2 + b_1(x_5^2+x_5x_6 + x_6^2) + b_2(x_7^2+x_7x_8 + x_8^2),
\end{equation}
where the coefficients $a_i \in \{1,2,3\}$, $1\le i\le 4$ and $b_1,b_2\in \{1,2,4\}$. 
Without loss of generality we can assume that $a_1\le a_2\le a_3\le a_4$ and $b_1\le b_2$. 
In \cite{a-a-l-w}, A. Alaca et. \!\!al  considered similar types of quadratic forms in four variables, which are either sums of four squares with coefficients $1,2,3,4$ or $6$ ($7$ such forms) or direct sum of the sums of two squares with coefficients $1$ or $3$ and the quadratic form $x^2+xy+y^2$ with coefficients $1, 2$ or $4$ ($6$ such forms). They used theta function identities to determine the representation formulas for these $13$  quadratic forms. 
In our recent work \cite{r-s-s}, we constructed  bases for the space of modular forms of weight $4$ for the group $\Gamma_0(48)$ with character, and used modular forms techniques to determine the number of representations of a natural number $n$ by certain octonary quadratic forms with coefficients $1,2,3,4,6$. Finding formulas for the number of representations for octonary quadratic forms with coefficients $1,2,3$ or $6$ were considered by various authors 
using several methods (see for example \cite{{a-a-w},{a-a-w1},{a-k-new},{a-k1},{a-k2},{a-w}}).  
In the present work, we adopt similar (modular forms) techniques to obtain the representation formulas. We show directly that the theta series corresponding to each of the quadratic form considered belongs to the space of modular forms of weight $4$ on $\Gamma_0(24)$ with some character (depending on the coefficients). Now, by constructing a basis for the space of modular forms $M_4(\Gamma_0(24),\chi)$  we find the required formulas. Here $\chi$ is either the trivial Dirichlet character modulo $24$ or one of the primitive Dirichlet characters (modulo $m$) $\chi_m = \left(\frac{m}{\cdot}\right)$, $m=8, 12, 24$. 
Since $M_4(\Gamma_0(24),\chi) \subseteq M_4(\Gamma_0(48),\chi)$, where $\chi$ is a Dirichlet character modulo $24$, we get the required explict bases from the basis of modular forms $M_4(\Gamma_0(48),\psi)$, where $\psi$ is a Dirichlet character modulo $48$, which was constructed in \cite{r-s-s}.

In the second part of the paper, we consider the quadratic forms of eight variables given by:
\begin{equation}\label{quad2}
(x_1^2+x_1x_2+x_2^2) + c_1(x_3^2+x_3x_4+x_4^2) + c_2(x_5^2+x_5x_6+x_6^2) + c_3(x_7^2+x_7x_8+x_8^2),
\end{equation}
where $c_1\le c_2\le c_3$ and $c_1,c_2,c_3\in \{1,2,4,8\}$. We note that for the $c_i$'s in the list, each of the quadratic form represents a theta series which belong to the space $M_4(\Gamma_0(24))$. Therefore, using our methods adopted for the earlier case, we also determine explicit formulas for the number of representations of a natural number by these class of quadratic forms. 

The total number of such quadratic forms given by \eqref{quad1} with coefficients $a_i\in \{1,2,3\}$ 
\linebreak 
and $b_i\in \{1,2,4\}$ is $90$. Each quadratic form in this list is denoted as 
a sextuple 
\linebreak 
$(a_1,a_2,a_3,a_4, b_1,b_2)$ and we list them in Table 1. We also put them in four classes corresponding to each of the modular forms space $M_4(\Gamma_0(24), \chi)$.   Similarly, we list the quadratic forms (total 19) given by \eqref{quad2} in Table 2. In this case all the corresponding theta series belong to $M_4(\Gamma_0(24))$. We do not consider the case 
$(1,1,1,1)$ as the formula is already known (see \cite[Theorem 17.4]{williams-book}).
It was shown that  $s_8(n) = 24 \sigma_3(n) + 216 \sigma_3(n/3)$.  
In our notation (see \S 3) $s_8(n) = M(1,1,1,1;n)$. 
Also, the cases $(1,2,2,4)$ and $(1,2,4,8)$ has been proved in \cite{kokluce} by using convolution sums method. 

The paper is organized as follows. In \S 2 we present the theorems proved in this article  and in \ \S 3 we give some preliminary results which are needed in proving the theorems. 
In \S 4 we give a proof of our theorems using the theory of modular forms. 

%\smallskip

\bigskip
\bigskip

\begin{center}
\textbf{Table 1.\\}
%\vskip 0.2 cm
List of quadratic forms in 8 variables given in \eqref{quad1}\\

\begin{tabular}{|c|c|}
\hline
$(a_1,a_2,a_3,a_4,b_1,b_2)$ & space \\
\hline
$(1,1,1,1,1,1),(1,1,1,1,1,2),(1,1,1,1,1,4),(1,1,1,1,2,2)$&~~\\
$(1,1,1,1,2,4),(1,1,1,1,4,4),(1,1,2,2,1,1),(1,1,2,2,1,2)$&~~\\
$(1,1,2,2,1,4),(1,1,2,2,2,2),(1,1,2,2,2,4),(1,1,2,2,4,4)$&~~\\
$(1,1,3,3,1,1),(1,1,3,3,1,2),(1,1,3,3,1,4),(1,1,3,3,2,2)$&~~\\
$(1,1,3,3,2,4),(1,1,3,3,4,4),(2,2,2,2,1,1),(2,2,2,2,1,2)$&$M_4(\Gamma_0(24))$\\
$(2,2,2,2,1,4),(2,2,2,2,2,2),(2,2,2,2,2,4),(2,2,2,2,4,4)$&~~\\
$(2,2,3,3,1,1),(2,2,3,3,1,2),(2,2,3,3,1,4),(2,2,3,3,2,2)$&~~\\
$(2,2,3,3,2,4),(2,2,3,3,4,4),(3,3,3,3,1,1),(3,3,3,3,1,2)$&~~\\
$(3,3,3,3,1,4),(3,3,3,3,2,2),(3,3,3,3,2,4),(3,3,3,3,4,4)$&~~\\
\hline
$(1,1,1,2,1,1),(1,1,1,2,1,2),(1,1,1,2,1,4),(1,1,1,2,2,2)$&~~\\
$(1,1,1,2,2,4),(1,1,1,2,4,4),(1,2,2,2,1,1,(1,2,2,2,1,2)$&~~\\
$(1,2,2,2,1,4),(1,2,2,2,2,2),(1,2,2,2,2,4),(1,2,2,2,4,4)$&$M_4(\Gamma_0(24),\chi_8)$\\
$(1,2,3,3,1,1),(1,2,3,3,1,2),(1,2,3,3,1,4),(1,2,3,3,2,2)$&~~\\
$(1,2,3,3,2,4),(1,2,3,3,4,4)$&~~\\
\hline
$(1,1,1,3,1,1),(1,1,1,3,1,2),(1,1,1,3,1,4),(1,1,1,3,2,2)$&~~\\
$(1,1,1,3,2,4),(1,1,1,3,4,4),(1,2,2,2,1,1,(1,2,2,3,1,2)$&~~\\
$(1,2,2,3,1,4),(1,2,2,3,2,2),(1,2,2,3,2,4),(1,2,2,3,4,4)$&$M_4(\Gamma_0(24),\chi_{12})$\\
$(1,3,3,3,1,1),(1,3,3,3,1,2),(1,3,3,3,1,4),(1,3,3,3,2,2)$&~~\\
$(1,3,3,3,2,4),(1,3,3,3,4,4)$&~~\\
\hline
$(1,1,2,3,1,1),(1,1,2,3,1,2),(1,1,2,3,1,4),(1,1,2,3,2,2)$&~~\\
$(1,1,2,3,2,4),(1,1,2,3,4,4),(2,2,2,3,1,1),(2,2,2,3,1,2)$&~~\\
$(2,2,2,3,1,4),(2,2,2,3,2,2),(2,2,2,3,2,4),(2,2,2,3,4,4)$&$M_4(\Gamma_0(24),\chi_{24})$\\
$(2,3,3,3,1,1),(2,3,3,3,1,2),(2,3,3,3,1,4),(2,3,3,3,2,2)$&~~\\
$(2,3,3,3,2,4),(2,3,3,3,4,4)$&~~\\
\hline
\end{tabular}
\end{center}

\bigskip

\newpage

\begin{center}
\textbf{Table 2.}
\vskip 0.2 cm
List of quadratic forms in \eqref{quad2} indicated by $(1,c_1,c_2,c_3)$.\\
%All the corresponding theta series lie in $M_4(\Gamma_0(24)$.
\begin{tabular}{|c|c|}
\hline
$(1,c_1,c_2,c_3)$& space\\
\hline
$(1,1,1,2),(1,1,1,4),(1,1,1,8),(1,1,2,2),(1,1,2,4),(1,1,2,8),(1,1,4,4)$&~~\\
$(1,1,4,8),(1,1,8,8),(1,2,2,2),(1,2,2,4),(1,2,2,8),(1,2,4,4),(1,2,4,8)$&$M_4(\Gamma_0(24))$\\
$(1,2,8,8),(1,4,4,4),(1,4,4,8),(1,4,8,8),(1,8,8,8)$&~~\\
\hline 
\end{tabular}
\end{center}

\section{Statement of results}

Let ${\mathbb N}, {\mathbb N}_0$ and ${\mathbb Z}$ denote the set of positive integers, non-negative integers and integers respectively. 
For $(a_1, a_2,a_3, a_4,b_1,b_2)$ as in Table 1, we define 
\begin{equation*}
\begin{split}
N(a_1,a_2,a_3,a_4,b_1,b_2;n) := \qquad \hskip 9.7cm &\\
\quad  \# \left\{(x_1,\ldots, x_8)\in {\mathbb Z}^8 \big\vert n = 
\sum_{i=1}^4a_ix_i^2 + b_1(x_5^2+x_5x_6 + x_6^2) + b_2(x_7^2+x_7x_8 + x_8^2)
\right\}.\\
\end{split}
\end{equation*}
to be the number of representations of $n$ by the quadratic form \eqref{quad1}.
Note that 
\linebreak 
$N(a_1,a_2,a_3,a_4,b_1,b_2;0) =1$. The formulas corresponding to Table 1 are stated in the following theorem. Formulas are divided into four parts each corresponding to one of the 
four spaces of modular forms $M_4(\Gamma_0(24),\chi)$. 

\begin{thm}\label{1}
Let $n\in {\mathbb N}$. \\
{\rm (i)} For each entry $(a_1,a_2,a_3,a_4,b_1,b_2)$ in Table {\rm 1} corresponding to the space $M_4(\Gamma_0(24))$, we have 
\begin{equation}
N(a_1,a_2,a_3,a_4,b_1,b_2;n)  = \sum_{i=1}^{16} \alpha_i A_i(n),
\end{equation}
where $A_i(n)$ are the Fourier coefficients of the basis elements $f_i$ defined in \S {\rm 4.1} and the values of the constants $\alpha_i$s are given in Table {\rm 3}. \\ 
{\rm (ii)} For each entry $(a_1,a_2,a_3,a_4,b_1,b_2)$ in Table {\rm 1} corresponding to the space $M_4(\Gamma_0(24),\chi_8)$, we have 
\begin{equation}
N(a_1,a_2,a_3,a_4,b_1,b_2;n)  = \sum_{i=1}^{14} \beta_i B_i(n),
\end{equation}
where $B_i(n)$ are the Fourier coefficients of the basis elements $g_i$ defined in \S {\rm 4.2} and the values of the constants $\beta_i$'s are given in Table {\rm 4}. \\
{\rm (iii)} For each entry $(a_1,a_2,a_3,a_4,b_1,b_2)$ in Table {\rm 1} corresponding to the space $M_4(\Gamma_0(24), \chi_{12})$, we have 
\begin{equation}
N(a_1,a_2,a_3,a_4,b_1,b_2;n)  = \sum_{i=1}^{16} \gamma_i C_i(n),
\end{equation}
where $C_i(n)$ are the Fourier coefficients of the basis elements $h_i$ defined in \S {\rm 4.3} and the values of the constants $\gamma_i$'s are given in Table {\rm 5}. \\
{\rm (iv)} For each entry $(a_1,a_2,a_3,a_4,b_1,b_2)$ in Table {\rm 1} corresponding to the space $M_4(\Gamma_0(24), \chi_{24})$, we have 
\begin{equation}
N(a_1,a_2,a_3,a_4,b_1,b_2;n)  = \sum_{i=1}^{14} \delta_i D_i(n),
\end{equation}
where $D_i(n)$ are the Fourier coefficients of the basis elements $F_i$ defined in \S {\rm 4.4} and the values of the constants $\delta_i$'s are given in Table {\rm 6}. 
\end{thm}

\bigskip

\noindent 
Now we consider the class of quadratic forms given by \eqref{quad2}. 
For $(1, c_1,c_2,c_3)$ as in Table 2, we define 
\begin{equation*}
\begin{split}
M(1, c_1,c_2,c_3;n) :=  \# \left\{(x_1,\ldots, x_8)\in {\mathbb Z}^8 \big\vert \right. \hskip 8cm &\\
\left. n =  \!(x_1^2+x_1x_2+x_2^2) + \!c_1(x_3^2+x_3x_4 + x_4^2) + \!c_2(x_5^2+x_5x_6
+ x_6^2) + \!c_3(x_7^2+x_7x_8+x_8^2)\right\}.&\\
\end{split}
\end{equation*}
to be the number of representations of $n$ by the quadratic form \eqref{quad2}.
Note that 
\linebreak 
$M(1,c_1,c_2,c_3;0) =1$. The formulas corresponding to Table 2 are stated in the following theorem. 

\begin{thm}\label{2}
Let $n\in {\mathbb N}$. \\
For each entry $(1,c_1, c_2,c_3;n)$ in Table {\rm 2}, we have 
\begin{equation}
M(1, c_1,c_2,c_3;n)  = \sum_{i=1}^{16} \nu_i A_i(n),
\end{equation}
where $A_i(n)$ are the Fourier coefficients of the basis elements $f_i$ defined in \S {\rm 4.1} and the values of the constants $\nu_i$s are given in Table {\rm 7}. \\ 
\end{thm}

\smallskip

\noindent 
{\bf Remark 2.1.}  Since one can write down the exact formulas using the explicit Fourier coefficients of the basis elements and using the coefficients tables given in each of the cases, we have not stated explicit formulas in the theorems (due to large number of such formulas). 
However, in \S 5 (at the end of the Tables), we give some sample formulas corresponding to each case. 

\smallskip

\section{Preliminaries}

In this section we present some preliminary facts on modular forms. For $k\in \frac{1}{2}{\mathbb Z}$, let $M_k(\Gamma_0(N),\chi)$ denote the space of modular forms of weight $k$ for the congruence subgroup $\Gamma_0(N)$ with character $\chi$ and $S_k(\Gamma_0(N), \chi)$ be the subspace of cusp forms of weight $k$ for $\Gamma_0(N)$ with character $\chi$. We assume $4\vert N$ when $k$ is not an integer and in that case, the character $\chi$ (which is a Dirichlet character modulo $N$) is an even character. When $\chi$ is the trivial (principal) character modulo $N$, we shall denote the spaces by $M_k(\Gamma_0(N))$ and $S_k(\Gamma_0(N))$ respectively. Further, when $k\ge 4$ is an integer and $N=1$, we shall denote these vector spaces by $M_k$ and $S_k$ respectively. 

For an integer $k \ge 4,$ let $E_k$ denote the normalized Eisenstein series of weight $k$ in $M_k$ given by 
$$
E_k(z) = 1 - \frac{2k}{B_k}\sum_{n\ge 1} \sigma_{k-1}(n) q^n,
$$
where $q=e^{2 i\pi z}$, $\sigma_r(n)$ is the sum of the $r$th powers of the positive divisors of $n$, and $B_k$ is the $k$-th Bernoulli number defined by $\displaystyle{\frac{x}{e^x-1} = \sum_{m=0}^\infty \frac{B_m}{m!} x^m}$.

The classical theta function which is fundamental to the theory of modular forms of half-integral weight is defined by 
\begin{equation}\label{theta}
\Theta(z) = \sum_{n\in {\mathbb Z}} q^{n^2},
\end{equation}
and is a modular form in the space $M_{1/2}(\Gamma_0(4))$. Another function which is mainly used in our work is the Dedekind eta function $\eta(z)$ and 
it is given by 
\begin{equation}\label{eta}
\eta(z)=q^{1/24} \prod_{n\ge1}(1-q^n).
\end{equation}

An eta-quotient is a finite product of integer powers of $\eta(z)$ and we denote it as follows: 
\begin{equation}\label{eta-q} 
\prod_{i=1}^s \eta^{r_i}(d_i z) := d_1^{r_1} d_2^{r_2} \cdots d_s^{r_s},
\end{equation}
where $d_i$'s are positive integers and $r_i$'s are non-zero integers.

We denote the theta series associated to the quadratic form $x^2+xy+y^2$ by 
\begin{equation}\label{F1}
{\mathcal F}(z) = \sum_{x,y\in {\mathbb Z}} q^{x^2+xy+y^2}.
\end{equation}
This function is referred to as the Borweins' two dimensional theta function in the literature.  
By \cite[Theorem 4]{schoeneberg}, it follows that ${\mathcal F}(z)$ is a modular form in 
$M_1(\Gamma_0(3),\chi_{-3})$. Here and in the sequel, for $m<0$, the character $\chi_m$ is the odd Dirichlet character modulo $|m|$ given by $\left(\frac{-m}{\cdot}\right)$. 

In the following we shall present some facts about modular forms of integral and half-integral weights, which we shall be using in our proof. We state them as lemmas, whose proofs follow from elementary theory of modular forms (of integral and half-integral weights).

\noindent {\bf  Lemma 1.} {\em (Duplication of modular forms)\\
%% We give a fact about certain duplication of modular forms. 
If $f$ is a modular form in $M_k(\Gamma_0(N), \chi)$, then for a positive integer $d$, the function $f(dz)$ is a modular form in $M_k(\Gamma_0(dN), \chi)$, 
if $k$ is an integer and it belongs to the space $M_k(\Gamma_0(dN), \chi \chi_d)$, if $k$ is a half-integer.
%%Here, we have used the notation $\chi_m := \left(\frac{m}{\cdot}\right)$, the Kronecker %%symbol, where $m$ is a non-zero integer, which is a 
% %character modulo $|m|$. 
}

\smallskip

\noindent {\bf Lemma 2.} {\em 
For  positive integers $r$, $r_1$, $r_2$,  $d_1$, $d_2$, we have 
\begin{equation}
\Theta^{r}(d_1z) \in 
\begin{cases} M_{r/2}(\Gamma_0(4d_1), \chi_{d_1}) & {\rm ~if~} r {\rm ~is~odd}, \\ 
                           M_{r/2}(\Gamma_0(4d_1), \chi_{-4}) & {\rm ~if~} r \equiv 2\pmod{4},\\ 
                           M_{r/2}(\Gamma_0(4d_1)) & {\rm ~if~} r \equiv 0\pmod{4}.
\end{cases}   
\end{equation}

\begin{equation}
\Theta^{r_1}(d_1z)\cdot \Theta^{r_2}(d_2z)  \in 
\begin{cases} 
M_{\frac{r_1+r_2}{2}}(\Gamma_0(4 [d_1, d_2]), \chi_{(-d_1d_2)}) & \!\!{\rm ~if~} r_1 r_2 {\rm ~is~odd~},  r_1+r_2\equiv 2\!\!\!\!\!\pmod{4}, \\
M_{\frac{r_1+r_2}{2}}(\Gamma_0(4 [d_1, d_2]), \chi_{(d_1d_2)}) &  \!\!{\rm ~if~} r_1 r_2 {\rm ~is~odd~}, r_1+r_2 \equiv 0\!\!\!\!\!\pmod{4}.
\end{cases}   
\end{equation}
}

\smallskip

\noindent {\bf Lemma 3.} {\em If $f_i \in M_{k_i}(\Gamma_0(M_i), \psi_i)$, $i=1,2$, then the 
product $f_1\cdot f_2$ is a modular form in $M_{k_1+k_2}(\Gamma_0(M), \psi_1\psi_2)$, where 
$M= {\rm lcm}(M_1,M_2)$. 
}

\smallskip

\noindent {\bf Lemma 4.} {\em 
The vector space $M_k(\Gamma_1(N))$ is decomposed as a direct sum:  
\begin{equation}
M_k(\Gamma_1(N)) = \oplus_{\chi}M_k(\Gamma_0(N), \chi),
\end{equation}
where the direct sum varies over all Dirichlet characters modulo $N$ if the weight $k$ is a  positive integer and varies over all even 
Dirichlet characters modulo $N$, $4\vert N$, if the weight $k$ is half-integer. Further, if $k$ is an integer, one has $M_k(\Gamma_0(N),\chi) = \{0\}$, if $\chi(-1) \not= (-1)^k$. We also have the following decomposition of the space into subspaces of Eisenstein series and cusp forms:
\begin{equation}
M_k(\Gamma_0(N),\chi) = {\mathcal E}_k(\Gamma_0(N),\chi) \oplus S_k(\Gamma_0(N), \chi),
\end{equation}
where ${\mathcal E}_k(\Gamma_0(N),\chi)$ is the space generated by the Eisenstein series of weight $k$ on $\Gamma_0(N)$ with character $\chi$. 
}

\smallskip

\noindent {\bf Lemma 5.} {\em 
By the Atkin-Lehner theory of newforms, the space $S_k(\Gamma_0(N),\chi)$ can be decomposed into the space of newforms and oldforms:
\begin{equation}
S_k(\Gamma_0(N),\chi) = S_k^{new}(\Gamma_0(N),\chi) \oplus S_k^{old}(\Gamma_0(N),\chi),
\end{equation}
where the above is an orthogonal direct sum (with respect to the Petersson scalar product) and  
\begin{equation}
S_k^{old}(\Gamma_0(N), \chi) = \bigoplus_{r\vert N, r<N\atop{rd\vert N}}S_k^{new}(\Gamma_0(r),\chi)\vert B(d).
\end{equation}
In the above, $S_k^{new}(\Gamma_0(N),\chi)$ is the space of newforms and $S_k^{old}(\Gamma_0(N),\chi)$ is the space of oldforms and the operator 
$B(d)$ is given by $f(z) \mapsto f(dz)$.
}

\smallskip

\noindent {\bf Lemma 6.} {\em 
Suppose that $\chi$ and $\psi$ are primitive Dirichlet characters with conductors $M$ and $N$, respectively. For a positive integer $k$, let 
\begin{equation}\label{eisenstein}
E_{k,\chi,\psi}(z) :=  c_0 + \sum_{n\ge 1}\left(\sum_{d\vert n} \psi(d) \cdot \chi(n/d) d^{k-1}\right) q^n,
\end{equation}
where 
$$
c_0 = \begin{cases}
0 &{\rm ~if~} M>1,\\
- \frac{B_{k,\psi}}{2k} & {\rm ~if~} M=1,
\end{cases}
$$
and $B_{k,\psi}$ denotes generalized Bernoulli number with respect to the character $\psi$. 
Then, the Eisenstein series $E_{k,\chi,\psi}(z)$ belongs to the space $M_k(\Gamma_0(MN), \chi/\psi)$, provided $\chi(-1)\psi(-1) = (-1)^k$ 
and $MN\not=1$. When $\chi=\psi =1$ (i.e., when $M=N=1$) and $k\ge 4$, we have $E_{k,\chi,\psi}(z) = - \frac{B_k}{2k} E_k(z)$, where $E_k$ is the normalized Eisenstein series of integer weight $k$ as defined before. We refer to \cite{{miyake}, {stein}} for details. 
}

\smallskip

We give a notation to the inner sum in \eqref{eisenstein}:
\begin{equation}\label{divisor}
\sigma_{k-1;\chi,\psi}(n) := \sum_{d\vert n} \psi(d) \cdot \chi(n/d) d^{k-1}.
\end{equation}

\smallskip

For more details on the theory of modular forms of integral and half-integral weights, we refer to \cite{{a-l}, {koblitz}, {li}, {miyake}, {schoeneberg}, {shimura}}.

\smallskip

\section{proofs of theorems}

In this section, we shall give a proof of our results. As mentioned in the introduction, we shall be using the theory of modular forms. 

The basic functions for the two types of quadratic forms considered in this paper are 
$\Theta(z)$ and ${\mathcal F}(z)$. To each quadratic form in \eqref{quad1} with coefficients $(a_1,a_2,a_3,a_4,b_1,b_2)$ as in Table 1, the associated theta series is given by 
\begin{equation*}
\Theta(a_1z) \Theta(a_2z) \Theta(a_3z) \Theta(a_4z) {\mathcal F}(b_1z) {\mathcal F}(b_2z).
\end{equation*} 
Using Lemma 1 and 2 along with the fact that ${\mathcal F}(z) \in 
M_1(\Gamma_0(3),\chi_{-3})$, it follows that the above product is a modular form in $M_4(\Gamma_0(24), \chi)$, 
where the character $\chi$ is one of the four characters that appear in Table 1 and it is determined by the coefficients $a_1,a_2,a_3,a_4$. As remarked earlier, the theta series corresponding to the form $x^2+xy+y^2$ is given by \eqref{F1} and it belongs to the 
space $M_1(\Gamma_0(3),\chi_{-3})$. Therefore, the associated modular form corresponding 
to the quadratic forms defined by \eqref{quad2} is given explicitly by 
\begin{equation*}
{\mathcal F}(z) {\mathcal F}(c_1z) {\mathcal F}(c_2z) {\mathcal F}(c_3z).
\end{equation*} 
Again by using Lemmas 1, 2 and 3 it follows that the above product is a modular form 
in $M_4(\Gamma_0(24))$. 
Therefore, in order to get the required formulae for $N(a_1,a_2,a_3,a_4,b_1,b_2;n)$ and 
$M(1, c_1,c_2,c_3;n)$  we need a basis for the above spaces of modular forms of level $24$. 
%%%%%
%We shall give explicit bases for the following spaces of modular forms of weight $4$:
%$$
%M_4(\Gamma_0(24)), M_4(\Gamma_0(24),\chi_8), M_4(\Gamma_0(24),\chi_{12}), %M_4(\Gamma_0(24),\chi_{24}).
%$$
(We have used the $L$-functions and modular forms database \cite{lfmdb} and \cite{martin} to get some of the cusp forms of weight $4$.)

\smallskip

%% \subsection{Proof of \thmref{1}}

\subsection{A basis for $M_4(\Gamma_0(24))$ and proof of \thmref{1}(i).}

The vector space $M_4(\Gamma_0(24))$ has dimension $16$ and we have $\dim_{\mathbb C}{\mathcal E}_4(\Gamma_0(24)) = 8$ and 
$\dim_{\mathbb C}S_4(\Gamma_0(24)) = 8$. For $d=6,8,12$ and $24$, $S_4^{new}(\Gamma_0(d))$ is one-dimensional. Let us define some eta-quotients and use them to give an explicit basis for $S_4(\Gamma_0(24))$.
Let\\
\begin{eqnarray}
f_{4,6}(z) = 1^{2}  2^{2} 3^{2}  6^2 := \displaystyle{\sum_{n\ge 1}} a_{4,6}(n) q^n, \quad 
%\eta^2(z) \eta^2(2z) \eta^2(3z) \eta^2(6z) 
f_{4,8}(z) = 2^4  4^4 := \displaystyle{\sum_{n\ge 1}} a_{4,8}(n) q^n, \\
%\eta^4(2z) \eta^4(4z) := 
f_{4,12}(z) =  
%\frac{\eta^2(2z)\eta^3(3z)\eta^3(4z)\eta^2(6z)}{\eta(z)\eta(12z)} - \frac{\eta^3(z) \eta^2(2z)\eta^2(6z)\eta^3(12z)}{\eta(3z)\eta(4z)} 
1^{-1} 2^2 3^3 4^3 6^2 12^{-1} - 1^3 2^2 3^{-1} 4^{-1} 6^2 12^3 := \displaystyle{\sum_{n\ge 1}} a_{4,12}(n) q^n, \\
%\frac{\eta^{16}(4z)}{\eta^4(2z)\eta^4(8z)} 
f_{4,24}(z) =  1^{-4} 2^{11} 3^{-4} 4^{-3} 6^{11} 12^{-3} := \displaystyle{\sum_{n\ge 1}} a_{4,24}(n) q^n.
%\frac{\eta^{11}(2z)\eta^{11}(6z)}{\eta^4(z)\eta^4(3z)\eta^3(4z)\eta^3(12z)} \\ 
\end{eqnarray}

We use the following notation in the sequel. For a Dirichlet character $\chi$ and a function $f$ with Fourier expansion $f(z) =\sum_{n\ge 1} a(n) q^n$, we define the twisted function $f \otimes \chi (z)$ as follows. 
\begin{equation}\label{twist}
f\otimes \chi (z) = \sum_{n\ge 1} \chi(n) a(n)q^n.
\end{equation}

A basis for the space $M_4(\Gamma_0(24))$ is given in the following proposition. 

\begin{prop}\label{trivial}
A basis for the Eisenstein series space ${\mathcal E}_4(\Gamma_0(24))$ is given by 
\begin{equation}\left\{E_4(tz), t\vert 24 \right\}
\end{equation}
and a basis for the space of cusp forms $S_4(\Gamma_0(24))$ is given by 
\begin{equation}
\begin{split}
\left\{f_{4,6}(t_1z), t_1\vert 4; f_{4,8}(t_2z), t_2\vert 3; f_{4,12}(t_3z), t_3\vert 2;f_{4,24}\otimes \chi_4(z)\right\}\\
\end{split}
\end{equation}
Together they form a basis for $M_4(\Gamma_0(24))$. 
\end{prop}
For the sake of simplicity in the formulae, we list these basis elements as $\{f_i(z)\vert 1\le i \le16\}$, where  $f_1(z) = E_4(z)$, $ f_2(z) = E_4(2z) $, $ f_3(z) = E_4(3z)$,
$ f_4(z) = E_4(4z)$, $f_5(z) = E_4(6z)$, $ f_6(z) = E_4(8z)$, $ f_7(z) = E_4(12z)$,$ f_{8}(z) = E_4(24z)$, $f_{9}(z) = f_{4,6}(z)$, $f_{10}(z) =  f_{4,6}(2z)$, $f_{11}(z) =  f_{4,6}(4z)$,  $f_{12}(z) = f_{4,8}(z)$, $f_{13}(z) = f_{4,8}(3z)$, $f_{14}(z) = f_{4,12}(z)$, $f_{15}(z) = f_{4,12}(2z)$, $f_{16}(z) =f_{4,24}\otimes \chi_4 (z)$\\

\noindent 
For $1\le i\le 16$, we denote the Fourier coefficients of the basis functions $f_i(z)$ as
$$
f_i(z) = \sum_{n\ge 1} A_i(n) q^n.
$$ 
We are now ready to prove the theorem. Noting that all the 36 cases corresponding to the trivial character in Table 1, the resulting functions belong to the space of modular forms of weight $4$ on $\Gamma_0(24)$ with trivial character (using Lemmas 1 to 3). So, we can express these theta functions as a linear combination of the basis given in \propref{trivial} as follows. 
\begin{equation}
\Theta(a_1z) \Theta(a_2z)\Theta(a_3z)\Theta(a_4z){\mathcal F}(b_1z) {\mathcal F}(b_2z) 
= \sum_{i=1}^{16} \alpha_i f_i(z),
\end{equation}
where $\alpha_i$'s are some explicit constants. Comparing the $n$-th Fourier coefficients on both the sides, we get 
\begin{equation*}
N(a_1,a_2,a_3,a_4, b_1,b_2;n) = \sum_{i=1}^{16} \alpha_i A_i(n).
\end{equation*}
Explicit values for the constants \!$\alpha_i$, \!$1\le i \le 16$ corresponding to these 36 cases are given in Table \!3. 

\smallskip

\subsection{A basis for $M_4(\Gamma_0(24),\chi_8)$ and proof of \thmref{1}(ii).}

The vector space $M_4(\Gamma_0(24),\chi_8)$ has dimension $14$ and we have $\dim_{\mathbb C}{\mathcal E}_4(\Gamma_0(24),\chi_8)) = 4$ and 
$\dim_{\mathbb C}S_4(\Gamma_0(24),\chi_8)) = 10$. For $d=6$ and $12$, $S_4^{new}(\Gamma_0(d),\chi_8) = \{0\}$. Also $S_4^{new}(\Gamma_0(8),\chi_8)$ is 2-dimensional and $S_4^{new}(\Gamma_0(24),\chi_8)$ is 6-dimensional.\\
In order to give explicit basis for this space,we define the following
\begin{equation}\label{eis:8chi8}
%\begin{split}
E_{4, {\bf 1},\chi_8}(z) ~=~ \frac{11}{2} + \sum_{n\ge 1} \sigma_{3;{\bf 1},\chi_8}(n) q^n, \quad E_{4,\chi_8, {\bf 1}}(z) ~=~  \sum_{n\ge 1} \sigma_{3;\chi_8, {\bf 1}}(n) q^n.
%\end{split}
\end{equation}
\begin{equation}\label{8chi8}
f_{4,8,\chi_8;1}(z) = 1^{-2}  2^{11} 4^{-3}  8^2 
%\frac{\eta^{11}(2z) \eta^2(8z)}{\eta^2(z) \eta^3(4z)}  
= \sum_{n\ge 1} a_{4,8,\chi_8;1}(n) q^n ,~~ 
f_{4,8,\chi_8;2}(z) = 1^{2}  2^{-3} 4^{11}  8^{-2}  
%\frac{\eta^{2}(z) \eta^{11}(4z)}{\eta^3(2z) \eta^2(8z)}
=  \sum_{n\ge 1} a_{4,8,\chi_8;2}(n) q^n. 
\end{equation}
For the space of Eisenstein series we use the basis elements of ${\mathcal E}_4(\Gamma_0(8),\chi_8)$ given in \eqref{eis:8chi8}.A basis for $S_4^{new}(\Gamma_0(8),\chi_8)$ is given in \eqref{8chi8}.  The following six eta-quotients span the space $S_4^{new}(\Gamma_0(24), \chi_8)$.
\begin{equation}
\begin{split}\label{24:chi8}
f_{4,24,\chi_8;1}(z) &=  1^2 2^1 3^{-4} 4^1 6^{10} 8^2 12^{-4} :=  \sum_{n\ge 1} a_{4,24,\chi_8;1}(n) q^n ,~~\\
%\frac{\eta^{2}(z) \eta(2z)\eta(4z) \eta^{10}(6z)\eta^{2}(8z)}{\eta^4(3z) \eta^4(12z)} 
f_{4,24,\chi_8;2}(z) & =  1^1 2^3 3^{-1} 4^1 6^4 8^{-1} 24^1 := \sum_{n\ge 1} a_{4,24,\chi_8;2}(n) q^n,~~\\
%\frac{\eta(z) \eta^3(2z)\eta(4z) \eta^{4}(6z)\eta(24z)}{\eta(3z) \eta(8z)}
f_{4,24,\chi_8;3}(z) & =  1^{-1} 2^4 3^1 6^3 8^1 12^1 24^{-1} := \sum_{n\ge 1} a_{4,24,\chi_8;3}(n) q^n , ~~\\
%\frac{\eta^{4}(2z) \eta(3z)\eta^{3}(6z) \eta(8z)\eta(12z)}{\eta(z) \eta(24z)} 
f_{4,24,\chi_8;4}(z) & = 1^{-2} 2^4 4^2 6^1 8^2 12^1  := \sum_{n\ge 1} a_{4,24,\chi_8;4}(n) q^n,\\
%\frac{\eta^{4}(2z) \eta^{2}(4z)\eta(6z) \eta^{2}(8z)\eta(12z)}{\eta^2(z)}
f_{4,24,\chi_8;5}(z) & = 2^1 3^{-2} 4^1 6^4 12^2 24^2 := \sum_{n\ge 1} a_{4,24,\chi_8;5}(n) q^n , ~~\\
%\frac{\eta(2z) \eta(4z)\eta^{4}(6z) \eta^{2}(12z)\eta^{2}(24z)}{\eta^2(3z)} 
f_{4,24,\chi_8;6}(z) & = 1^{-6} 2^{14} 6^1 8^{-2} 12^1  := \sum_{n\ge 1} a_{4,24,\chi_8;6}(n) q^n
%\frac{\eta^{14}(2z) \eta(6z)\eta(12z)}{\eta^6(z) \eta^2(8z)}.
\end{split}
\end{equation}
A basis for the space $M_4(\Gamma_0(24),\chi_{8})$ is given in the following proposition. 

\begin{prop}\label{chi8} 
A basis for the space $M_4(\Gamma_0(24),\chi_8)$ is given by 
\begin{equation}
\begin{split}
&\left\{E_{4, {\bf 1},\chi_8}(tz),~ E_{4,\chi_8, {\bf 1}}(tz), t\vert3; f_{4,8,\chi_8;1}(t_1z),f_{4,8,\chi_8;2}(t_1z), t_1\vert 3;  f_{4,24,\chi_8;1}(z),\right.\\
& \left.f_{4,24,\chi_8;2}(z),f_{4,24,\chi_8;3}(z),f_{4,24,\chi_8;4}(z),  f_{4,24,\chi_8;5}(z),f_{4,24,\chi_8;6}(z)\right\}\\
\end{split}
\end{equation}
where $E_{4,{\bf 1},\chi_8}(z)$ and $E_{4,\chi_8,{\bf 1}}(z)$ are defined in \eqref{eis:8chi8}, $f_{4,8,\chi_8;i}(z)$, $i=1,2$ are defined in \eqref{8chi8} and 
$f_{4,24,\chi_8;j}(z)$, $1\le j\le 6$ are defined by \eqref{24:chi8}.
\end{prop}
For the sake of simplifying the notation, we shall list the basis in \propref{chi8} as
$$
g_i(z) = \sum_{n\ge 1} B_i(n) q^n,  ~1\le i \le 14,
$$
where 
$g_1(z) = E_{4, {\bf 1},\chi_8}(z)$, $g_2(z) = E_{4, {\bf 1},\chi_8}(3z)$,
$g_3(z)=E_{4,\chi_8, {\bf 1}}(z)$,  $g_4(z) = E_{4,\chi_8, {\bf 1}}(3z)$,
$g_5(z) = f_{4,8,\chi_8;1}(z)$, $g_{6}(z) = f_{4,8,\chi_8;1}(3z)$, $g_{7}(z) = f_{4,8,\chi_8;2}(z)$,$g_{8}(z) = f_{4,8,\chi_8;2}(3z)$,$g_{9}(z)=$ 
%\break 
$f_{4,24,\chi_8;1}(z)$,
$g_{10}(z) = f_{4,24,\chi_8;2}(z)$, 
$g_{11}(z) = f_{4,24,\chi_8;3}(z)$, 
$g_{12}(z) = f_{4,24,\chi_8;4}(z)$,  
$g_{13}(z) = f_{4,24,\chi_8;5}(z)$, 
$g_{14}(z) = f_{4,24,\chi_8;6}(z)$,

\noindent 
We now prove \thmref{1}(ii). In this case, for all the 18 sextuples corresponding to the $\chi_8$ character space (in Table 1), the resulting products of theta functions are modular forms of weight $4$ on $\Gamma_0(24)$ with character $\chi_8$ (By Lemma 1 to 3). So, we can express these products of theta functions as a linear combination of the basis given in \propref{chi8}:  
\begin{equation}
\Theta(a_1z) \Theta(a_2z)\Theta(a_3z)\Theta(a_4z) {\mathcal F}(b_1z) {\mathcal F}(b_2z) = \sum_{i=1}^{14} \beta_i g_i(z).
\end{equation}
Comparing the $n$-th Fourier coefficients on both the sides, we get 
\begin{equation*}
N(a_1,a_2,a_3,a_4, b_1,b_2;n) = \sum_{i=1}^{14} \beta_i B_i(n).
\end{equation*}
Explicit values for the constants $\beta_i$, $1\le i \le 14$ corresponding to these 18 cases  
are given in Table 4.

\smallskip

\subsection{A basis for $M_4(\Gamma_0(24), \chi_{12})$ and proof of \thmref{1}(iii).} 

The dimension of the space in this case is $16$, with $\dim_{\mathbb C}{\mathcal E}_4(\Gamma_0(24), \chi_{12}) = 8$ and 
$\dim_{\mathbb C}S_4(\Gamma_0(24), \chi_{12}) = 8$. The old class is spanned by 
the space $S_4^{new}(\Gamma_0(12), \chi_{12})$, which is $4$ dimensional with spanning functions given by  the following four eta-quotients:
\begin{equation}
\begin{split}
f_{4,12,\chi_{12};1}(z) ~=~ 2^{-1} 3^4 4^2 6^5 12^{-2}, & \quad  
%% \frac{\eta^{4}(3z) \eta^{2}(4z)\eta^{5}(6z)}{\eta(2z) \eta^2(12z)} = \sum_{n\ge 1} a_{4,12,\chi_{12};1}(n) q^n ,&~~\\
f_{4,12,\chi_{12};2}(z) ~=~ 
%% \frac{\eta^4(3z)\eta^3(4z) \eta^{3}(12z)}{\eta^2(6z)}= \sum_{n\ge 1} a_{4,12,\chi_{12};2}(n) q^n,~~\\
3^4 4^3 6^{-2} 12^3, \\
f_{4,12,\chi_{12};3}(z) ~=~ 2^2 3^4 4^{-1} 6^{-4} 12^7, & \quad 
%% \frac{\eta^{2}(2z) \eta^4(3z)\eta^{7}(12z)}{\eta(4z)\eta^4(6z)} = \sum_{n\ge 1} a_{4,12,\chi_{12};3}(n) q^n , & ~~\\
f_{4,12,\chi_{12};4}(z) = 
%% \frac{\eta^{4}(z) \eta^{7}(12z)}{\eta(4z)\eta^2(6z)}= \sum_{n\ge 1} a_{4,12,\chi_{12};4}(n) q^n,\\
1^4 4^{-1} 6^{-2} 12^7. \\
\end{split}
\end{equation}
We write the Fourier expansions of these forms as $f_{4,12,\chi_{12};j}(z) = \displaystyle{\sum_{n\ge 1}} a_{4,12,\chi_{12};j}(n) q^n$, $1\le j\le 4$.\\
In the following proposition we give a basis for the space $M_4(\Gamma_0(24),\chi_{12})$. 

\begin{prop}\label{chi12}
A basis for the space $M_4(\Gamma_0(24),\chi_{12})$ is given by 
\begin{equation}
\begin{split}
\left\{E_{4, {\bf 1},\chi_{12}}(tz),{E_{4,\chi_{12}, {\bf 1}}(tz),E_{4,\chi_{-4},\chi_{-3}}(tz),E_{4,\chi_{-3},\chi_{-4}}(tz)}, t\vert 2; 
f_{4,12,\chi_{12};j}(t_1z), t_1\vert 2, 1\le j\le 4\right\},&\\
%f_{4,12,\chi_{12};2}(t_1z), f_{4,12,\chi_{12};3}(t_1z), f_{4,12,%\chi_{12};4}(t_1z),t_1\vert 4},&\\
\end{split}
\end{equation}
where the Eisenstein series in the basis are defined by \eqref{eisenstein}.
\end{prop}
Let us denote the 16 basis elements in the above proposition as follows. \\
$\left\{h_i(z)\vert 1\le i \le 16\right\}$, where $h_1(z) = E_{4, {\bf 1},\chi_{12}}(z)$, $h_2(z) = E_{4,\chi_{12}, {\bf 1}}(z)$, 
$h_3(z) = E_{4,\chi_{-4},\chi_{-3}}(z)$, $h_4(z) = E_{4,\chi_{-3},\chi_{-4}}(z)$,  $h_5(z) = E_{4, {\bf 1},\chi_{12}}(2z)$, $h_6(z) = E_{4,\chi_{12}, {\bf 1}}(2z)$,$h_7(z) = E_{4,\chi_{-4},\chi_{-3}}(2z)$, $h_8(z) = E_{4,\chi_{-3},\chi_{-4}}(2z)$,$h_{8+j}(z) =  f_{4,12,\chi_{12};j}(z)$, 
$1\le j\le 4$, $h_{12+j}(z) =  f_{4,12,\chi_{12};j}(2z)$, $1\le j\le 4$. 

\smallskip

To prove \thmref{1}(iii), we consider the case of  18 sextuples corresponding to the character $\chi_{12}$ in Table 1. The resulting products of theta functions are modular forms of weight $4$ on $\Gamma_0(24)$ with character $\chi_{12}$ (once again we  use Lemams 1 to 3 to get this). So, we can express each of these products of theta functions as a linear combination of the basis given in \propref{chi12} as follows. 
\begin{equation}
\Theta(a_1z) \Theta(a_2z)\Theta(a_3z)\Theta(a_4z) {\mathcal F}(b_1z) {\mathcal F}(b_2z) = \sum_{i=1}^{16} \gamma_i g_i(z).
\end{equation}
Comparing the $n$-th Fourier coefficients on both the sides, we get 
\begin{equation*}
N(a_1,a_2,a_3,a_4, b_1,b_2;n) = \sum_{i=1}^{16} \gamma_i C_i(n).
\end{equation*}
Explicit values of the constants $\gamma_i$, $1\le i \le 16$ corresponding to these 18 cases   are given in Table 5. 

\smallskip

%\smallskip

\subsection{A basis for $M_4(\Gamma_0(24), \chi_{24})$ and proof of \thmref{1}(iv).} ~ 
We have 
\linebreak 
$\dim_{\mathbb C} M_4(\Gamma_0(24), \chi_{24}) = 14$ and $\dim_{\mathbb C}{\mathcal E}_4(\Gamma_0(24),\chi_{24}) = 4$. To get the span of the Eisenstein series space 
${\mathcal E}_4(\Gamma_0(24),\chi_{24})$, we use the Eisenstein series $E_{4,\chi,\psi}(z)$ defined in \eqref{eisenstein}, where $\chi,\psi \in \{{\bf 1}, \chi_{-8}, \chi_{-12}, \chi_{24}\}$. Note that  for d= 6,8 and 12, $S_4^{new}(\Gamma_0(d),\chi_{24}) = \{0\}$ and the space $S_4^{new}(\Gamma_0(24),\chi_{24})$ is spanned by the following ten eta-quotients (notation as in \eqref{eta-q}):
\begin{equation}
\begin{split}
 f_{4,24,\chi_{24};1}(z) ~=~ 3^{-2} 6^7 8^3 12^3 24^{-3}, &\quad 
%% \frac{\eta^{7}(6z) \eta^{3}(8z)\eta^{3}(12z)}{\eta^2(3z) \eta^3(24z)} = \sum_{n\ge 1} a_{4,24,\chi_{24};1}(n) q^n ,&~~\\
f_{4,24,\chi_{24};2}(z) ~=~ 3^{2} 4^7 6^{-3} 8^{-2} 12^{4}, \\
%% \frac{\eta^2(3z)\eta^7(4z) \eta^{4}(12z)}{\eta^3(6z) \eta^2(8z)}= \sum_{n\ge 1} a_{4,24,\chi_{24};2}(n) q^n,~~\\
~~f_{4,24,\chi_{24};3}(z) = 3^{2} 4^{-3} 6^{1} 8^{6} 12^{2}, \qquad &\quad
%%\frac{\eta^{2}(3z) \eta(6z)\eta^{6}(8z)\eta^2(12z)}{\eta^3(4z)} = \sum_{n\ge 1} a_{4,24,\chi_{24};3}(n) q^n , & ~~\\
 f_{4,24,\chi_{24};4}(z) =  3^2 6^{-3} 8^3 12^5 24^1, \\
%%\frac{\eta^{2}(3z) \eta^{3}(8z)\eta^5(12z) \eta(24z)}{\eta^3(6z)}= \sum_{n\ge 1} a_{4,24,\chi_{24};4}(n) q^n,\\
f_{4,24,\chi_{24};5}(z) = 3^2 4^2 6^{-3} 8^{-1} 12^3 24^5, & \quad
%%\frac{\eta^2(3z) \eta(6z)\eta^{2}(8z) \eta^{4}(24z)}{\eta(4z)} = \sum_{n\ge 1} a_{4,24,\chi_{24};5}(n) q^n ,& ~~\\
 f_{4,24,\chi_{24};6}(z) = 3^2 4^1 6^1 8^{-2} 12^{-2} 24^8, \\ 
%%\frac{\eta^{2}(3z) \eta^2(4z)\eta^3(12z)\eta^5(24z)}{\eta^3(6z)\eta(8z)}= \sum_{n\ge 1} a_{4,24,\chi_{24};6}(n) q^n ,& ~~\\
 f_{4,24,\chi_{24};7}(z) = 3^2 4^1 6^1 8^{-2} 12^{-2} 24^8, & \quad
%%\frac{\eta^{2}(3z) \eta(4z)\eta(6z) \eta^{8}(24z)}{\eta^2(8z) \eta^2(12z)} = \sum_{n\ge 1} a_{4,24,\chi_{24};7}(n) q^n ,&~~\\
 f_{4,24,\chi_{24};8}(z) = 1^1 3^{-1} 6^1 8^{-2} 12^1 24^8, \\
%%\frac{\eta(z) \eta(6z)\eta(12z) \eta^{8}(24z)}{\eta(3z) \eta^2(8z)}= \sum_{n\ge 1} a_{4,24,\chi_{24};8}(n) q^n,~~\\
 f_{4,24,\chi_{24};9}(z) = 2^2 3^6 4^1 6^{-3} 8^2, &\quad 
%%\frac{\eta^{2}(2z) \eta^6(3z)\eta(4z)\eta^2(8z)}{\eta^3(6z)} = \sum_{n\ge 1} a_{4,24,\chi_{24};9}(n) q^n , & ~~\\
f_{4,24,\chi_{24};10}(z) = 3^2 4^3 6^5 12^{-4} 24^2. \\
%\frac{\eta^{2}(3z) \eta^{3}(4z)\eta^5(6z)\eta^{2}(24z)}{\eta^4(12z)}= \sum_{n\ge 1} a_{4,24,\chi_{24};10}(n) q^n,\\
\end{split}
\end{equation}
We write the Fourier expansions as $f_{4,24,\chi_{24};j}(z) = \sum_{n \ge 1} a_{4,24,\chi_{24};j}(n) q^n$.  We now give a basis for the space $M_4(\Gamma_0(24),\chi_{24})$ in the following proposition.

\begin{prop}\label{chi24}
The following functions span the space $M_4(\Gamma_0(24),\chi_{24})$. 
\begin{equation}
\begin{split}
& \quad \left\{E_{4, {\bf 1},\chi_{24}}(z),{E_{4,\chi_{24}, {\bf 1}}(z),E_{4,\chi_{-8},\chi_{-3}}(z),E_{4,\chi_{-3},\chi_{-8},}(z)},f_{4,24,\chi_{24};j}(z), 1\le j\le 10;\right\}.  
%f_{4,24,\chi_{24};3}(tz),f_{4,24,\chi_{24};4}(tz), f_{4,24,\chi_{24};5}(tz),\right. &\\
%\left.f_{4,24,\chi_{24};6}(tz),f_{4,24,\chi_{24};7}(tz),f_{4,24,\chi_{24};8}(tz),f_{4,24,\chi_{24};9}(tz),f_{4,24,\chi_{24};10}(tz); t\vert 2\right\}&\\
\end{split}
\end{equation}
\end{prop}

We list these basis elements as $\{F_i(z)\vert 1\le i \le 14\}$, where  
$ F_1(z) = E_{4, {\bf 1},\chi_{24}}(z)$,
$ F_2(z) = E_{4,\chi_{24}, {\bf 1}}(z)$,
$ F_3(z) = E_{4,\chi_{-8},\chi_{-3}}(z)$,
$ F_4(z) = E_{4,\chi_{-3},\chi_{-8},}(z)$,
$ F_5(z) = f_{4,24,\chi_{24};1}(z)$,
$ F_6(z) = f_{4,24,\chi_{24};2}(z)$, 
$ F_7(z) = f_{4,24,\chi_{24};3}(z)$,
$ F_8(z) = f_{4,24,\chi_{24};4}(z)$, 
$ F_9(z) = f_{4,24,\chi_{24};5}(z)$,
$ F_{10}(z) = f_{4,24,\chi_{24};6}(z)$, 
$ F_{11}(z) = f_{4,24,\chi_{24};7}(z)$,
$ F_{12}(z) =f_{4,24,\chi_{24};8}(z)$,
$ F_{13}(z) = f_{4,24,\chi_{24};9}(z)$, 
$ F_{14}(z) = f_{4,24,\chi_{24};10}(z)$,\\
\noindent 
As in the previous cases, we denote the Fourier coefficients of these basis functions by 
$$
F_i(z) = \sum_{n\ge 1} D_i(n) q^n, ~1\le i \le 14.
$$ 
%% \noindent 
To get the formula in \thmref{1}(iv), we note that for all the 18 sextuples corresponding to the character $\chi_{24}$ in Table 1, the resulting functions belong to the space 
$M_4(\Gamma_0(24), \chi_{24})$, by using Lemmas 1 to 3.  So, as before, we express these theta functions as linear combinations of the basis elements:
\begin{equation}
\Theta(a_1z) \Theta(a_2z)\Theta(a_3z)\Theta(a_4z) {\mathcal F}(b_1z) {\mathcal F}(b_2z) = \sum_{i=1}^{14} \delta_i g_i(z).
\end{equation}
Comparing the $n$-th Fourier coefficients on both the sides, we get 
\begin{equation*}
N(a_1,a_2,a_3,a_4, b_1,b_2;n) = \sum_{i=1}^{14} \delta_i D_i(n).
\end{equation*}
Explicit values of the constants $\delta_i$, $1\le i \le 14$ corresponding to these 18 cases  
corresponding to character $\chi_{24}$ are given in Table 6.

\smallskip

\subsection{Proof of \thmref{2}}
This theorem is corresponding to Table 2 and in this case all the product functions 
$$
{\mathcal F}(z) {\mathcal F}(c_1z) {\mathcal F}(c_2z) {\mathcal F}(c_3z)
$$
belong to the space $M_4(\Gamma_0(24))$. Therefore, proceeding as in the proof of 
\thmref{1}(i),  we express these theta functions as  linear combinations of the basis elements:
\begin{equation} 
{\mathcal F}(z) {\mathcal F}(c_1z) {\mathcal F}(c_2z) {\mathcal F}(c_3z) = \sum_{i=1}^{16} \nu_i f_i(z).
\end{equation}
Comparing the $n$-th Fourier coefficients on both the sides, we get 
\begin{equation*}
M(1,c_1,c_2,c_3;n) = \sum_{i=1}^{16} \nu_i A_i(n).
\end{equation*}
The constants $\nu_i$, $1\le i \le 16$ corresponding to the 19 cases  of table 2 are given in Table 7.

%\smallskip

\bigskip

\section{List of tables and sample formulas }

In this section we list the remaining tables mentioned in the theorems and provide explicit 
sample formulas in some cases. In the first subsection we list the tables and in the second subsection we give the sample formulas. 

\smallskip

\subsection{List of tables}
In this section, we list the tables 3, 4, 5, 6 and 7 which give the explicit coefficients that appear in the formulas of \thmref{1} and \thmref{2}. 

\bigskip

\newpage

\begin{center}
{\tiny 
\textbf{Table 3.} (Theorem 2.1 (i))
\begin{tabular}{|c|c|c|c|c|c|c|c|c|c|c|c|c|c|c|c|c|}
\hline
$(a_1a_2a_3a_4,b_1b_2)$&$ \alpha_1$ &$ \alpha_2$& $ \alpha_3$& $ \alpha_4$&$ \alpha_5$ &$ \alpha_6$& $ \alpha_7$& $ \alpha_8$&$ \alpha_9$ &$ \alpha_{10}$ & $ \alpha_{11}$&$ \alpha_{12}$& $ \alpha_{13}$&$ \alpha_{14}$& $ \alpha_{15}$&$ \alpha_{16}$\\
\hline 
&&&&&&&&&&&&&&&&\\
$(1111,11)$&$\frac{7}{75}$&$\frac{-7}{100}$&$\frac{-9}{25}$&$\frac{-28}{75}$& $\frac{27}{100}$&0&$\frac{36}{25}$&0&$\frac{-72}{5}$&$\frac{-288}{5}$&0&0&0&12&0&0\\
&&&&&&&&&&&&&&&&\\
$(1111,12)$
&$\frac{13}{300}$&$\frac{-13}{200}$&$\frac{9}{100}$&$\frac{26}{75}$&$\frac{-27}{200}$&0&$\frac{18}{25}$&0&$\frac{48}{5}$&$\frac{96}{5}$&0&0&0&-6&0&0\\
&&&&&&&&&&&&&&&&\\
$(1111,14)$
&$\frac{7}{300}$&0&$\frac{-9}{100}$&$\frac{-28}{75}$&0&0&$\frac{36}{25}$&0&$\frac{-18}{5}$&$\frac{72}{5}$&0&0&0&12&0&0\\
&&&&&&&&&&&&&&&&\\
$(1111,22)$
&$\frac{7}{300}$&0&$\frac{-9}{100}$&$\frac{-28}{75}$&0&0&$\frac{36}{25}$&0&$\frac{12}{5}$&$\frac{-48}{5}$&0&0&0&0&0&0\\
&&&&&&&&&&&&&&&&\\
$(1111,24)$
&$\frac{13}{1200}$&$\frac{-13}{400}$&$\frac{9}{400}$&$\frac{26}{75}$&$\frac{-27}{400}$&0&$\frac{18}{25}$&0&$\frac{12}{5}$&$\frac{96}{5}$&0&0&0&3&0&0\\
&&&&&&&&&&&&&&&&\\
$(1111,44)$ 
&$\frac{7}{1200}$&$\frac{7}{400}$&$\frac{-9}{400}$&$\frac{-28}{75}$&$\frac{-27}{400}$&0&$\frac{36}{25}$&0&$\frac{18}{5}$&$\frac{72}{5}$&0&0&0&3&0&0\\
&&&&&&&&&&&&&&&&\\
$(1122,11)$ 
&$\frac{7}{150}$&$\frac{-7}{150}$&$\frac{-9}{50}$&$\frac{7}{300}$&$\frac{9}{50}$&$\frac{-28}{75}$&$\frac{-9}{100}$&$\frac{36}{25}$&$\frac{-36}{5}$&-48&$\frac{-768}{5}$&-3&-81&6&36&9\\
&&&&&&&&&&&&&&&&\\
$(1122,12)$ 
&$\frac{13}{600}$&$\frac{-13}{600}$&$\frac{9}{200}$&$\frac{-13}{600}$&$\frac{-9}{200}$&$\frac{26}{75}$&$\frac{-9}{200}$&$\frac{18}{25}$&$\frac{24}{5}$&12&$\frac{96}{5}$&$\frac{15}{2}$&$\frac{81}{2}$&-3&-6&$\frac{-9}{2}$\\
&&&&&&&&&&&&&&&&\\
$(1122,14)$ 
&$\frac{7}{600}$&$\frac{-7}{600}$&$\frac{-9}{200}$&$\frac{7}{300}$&$\frac{9}{200}$&$\frac{-28}{75}$&$\frac{-9}{100}$&$\frac{36}{25}$&$\frac{-9}{5}$&6&$\frac{-48}{5}$&$\frac{-3}{2}$&$\frac{-81}{2}$&6&0&$\frac{9}{2}$\\
&&&&&&&&&&&&&&&&\\
$(1122,22)$
&$\frac{7}{600}$&$\frac{-7}{600}$&$\frac{-9}{200}$&$\frac{7}{300}$&$\frac{9}{200}$&$\frac{-28}{75}$&$\frac{-9}{100}$&$\frac{36}{25}$&$\frac{6}{5}$&-12&$\frac{-288}{5}$&0&0&0&12&0\\
&&&&&&&&&&&&&&&&\\
$(1122,24)$ 
&$\frac{13}{2400}$&$\frac{-13}{2400}$&$\frac{9}{800}$&$\frac{-13}{600}$&$\frac{-9}{800}$&$\frac{26}{75}$&$\frac{-9}{200}$&$\frac{18}{25}$&$\frac{6}{5}$&12&$\frac{96}{5}$&0&0&$\frac{3}{2}$&-6&0\\
&&&&&&&&&&&&&&&&\\
$(1122,44)$ 
&$\frac{7}{2400}$&$\frac{-7}{2400}$&$\frac{-9}{800}$&$\frac{7}{300}$&$\frac{9}{800}$&$\frac{-28}{75}$&$\frac{-9}{100}$&$\frac{36}{25}$&$\frac{9}{5}$&6&$\frac{-48}{5}$&0&0&$\frac{3}{2}$&0&0\\
&&&&&&&&&&&&&&&&\\
$(1133,11)$ 
&$\frac{2}{75}$&$\frac{-1}{30}$&$\frac{6}{25}$&$\frac{8}{75}$&$\frac{-3}{10}$&0&$\frac{24}{25}$&0&$\frac{48}{5}$&$\frac{288}{5}$&0&0&0&0&0&0\\
&&&&&&&&&&&&&&&&\\
$(1133,12)$ 
&$\frac{1}{60}$&$\frac{-1}{120}$&$\frac{-3}{20}$&$\frac{-2}{15}$&$\frac{3}{40}$&0&$\frac{6}{5}$&0&0&0&0&0&0&6&0&0\\
&&&&&&&&&&&&&&&&\\
$(1133,14)$ 
&$\frac{1}{150}$&$\frac{-1}{75}$&$\frac{3}{50}$&$\frac{8}{75}$&$\frac{-3}{25}$&0&$\frac{24}{25}$&0&$\frac{42}{5}$&$\frac{168}{5}$&0&0&0&0&0&0\\
&&&&&&&&&&&&&&&&\\
$(1133,22)$ 
&$\frac{1}{150}$&$\frac{-1}{75}$&$\frac{3}{50}$&$\frac{8}{75}$&$\frac{-3}{25}$&0&$\frac{24}{25}$&0&$\frac{12}{5}$&$\frac{48}{5}$&0&0&0&0&0&0\\
&&&&&&&&&&&&&&&&\\
$(1133,24)$
&$\frac{1}{240}$&$\frac{1}{240}$&$\frac{-3}{80}$&$\frac{-2}{15}$&$\frac{-3}{80}$&0&$\frac{6}{5}$&0&0&0&0&0&0&3&0&0\\
&&&&&&&&&&&&&&&&\\
$(1133,44)$
&$\frac{1}{600}$&$\frac{-1}{120}$&$\frac{3}{200}$&$\frac{8}{75}$&$\frac{-3}{40}$&0&$\frac{24}{25}$&0&$\frac{18}{5}$&$\frac{48}{5}$&0&0&0&0&0&0\\
&&&&&&&&&&&&&&&&\\
$(2222,11)$ 
&$\frac{7}{400}$&$\frac{91}{1200}$&$\frac{-27}{400}$&$\frac{-7}{100}$&$\frac{-117}{400}$&$\frac{-28}{75}$&$\frac{27}{100}$&$\frac{36}{25}$&$\frac{-36}{5}$&$\frac{-312}{5}$&$\frac{-768}{5}$&-3&-81&9&36&9\\
&&&&&&&&&&&&&&&&\\
$(2222,12)$ 
&$\frac{13}{800}$&$\frac{-247}{2400}$&$\frac{27}{800}$&$\frac{13}{200}$&$\frac{-171}{800}$&$\frac{26}{75}$&$\frac{27}{200}$&$\frac{18}{25}$&$\frac{18}{5}$&$\frac{84}{5}$&$\frac{96}{5}$&$\frac{15}{2}$&$\frac{81}{2}$&$\frac{-9}{2}$&-6&$\frac{-9}{2}$\\
&&&&&&&&&&&&&&&&\\
$(2222,14)$ 
&$\frac{7}{800}$&$\frac{-133}{2400}$&$\frac{-27}{800}$&$\frac{7}{100}$&$\frac{171}{800}$&$\frac{-28}{75}$&$\frac{-27}{100}$&$\frac{36}{25}$&$\frac{-18}{5}$&$\frac{-24}{5}$&$\frac{-48}{5}$&$\frac{-3}{2}$&$\frac{-81}{2}$&$\frac{9}{2}$&0&$\frac{9}{2}$\\
&&&&&&&&&&&&&&&&\\
$(2222,22)$ 
&0&$\frac{7}{75}$&0&$\frac{-7}{100}$&$\frac{-9}{25}$&$\frac{-28}{75}$&$\frac{27}{100}$&$\frac{36}{25}$&0&$\frac{-72}{5}$&$\frac{-288}{5}$&0&0&0&12&0\\
&&&&&&&&&&&&&&&&\\
$(2222,24)$ 
&0&$\frac{13}{300}$&0&$\frac{-13}{200}$&$\frac{9}{100}$&$\frac{26}{75}$&$\frac{-27}{200}$&$\frac{18}{25}$&0&$\frac{48}{5}$&$\frac{96}{5}$&0&0&0&-6&0\\
&&&&&&&&&&&&&&&&\\
$(2222,44)$ 
&0&$\frac{7}{300}$&0&0&$\frac{-9}{100}$&$\frac{-28}{75}$&0&$\frac{36}{25}$&0&$\frac{12}{5}$&$\frac{-48}{5}$&0&0&0&0&0\\
&&&&&&&&&&&&&&&&\\
$(2233,11)$ 
&$\frac{1}{75}$&$\frac{-1}{75}$&$\frac{3}{25}$&$\frac{-1}{150}$&$\frac{-3}{25}$&$\frac{8}{75}$&$\frac{-3}{50}$&$\frac{24}{25}$&$\frac{24}{5}$&48&$\frac{768}{5}$&10&18&0&-24&-6\\
&&&&&&&&&&&&&&&&\\
$(2233,12)$
&$\frac{1}{120}$&$\frac{-1}{120}$&$\frac{-3}{40}$&$\frac{1}{120}$&$\frac{3}{40}$&$\frac{-2}{15}$&$\frac{-3}{40}$&$\frac{6}{5}$&0&-12&-96&$\frac{-7}{2}$&$\frac{-9}{2}$&3&6&$\frac{9}{2}$\\
&&&&&&&&&&&&&&&&\\
$(2233,14)$ 
&$\frac{1}{300}$&$\frac{-1}{300}$&$\frac{3}{100}$&$\frac{-1}{150}$&$\frac{-3}{100}$&$\frac{8}{75}$&$\frac{-3}{50}$&$\frac{24}{25}$&$\frac{21}{5}$&18&$\frac{48}{5}$&$\frac{5}{2}$&$\frac{63}{2}$&0&-12&$\frac{-3}{2}$\\
&&&&&&&&&&&&&&&&\\
$(2233,22)$ 
&$\frac{1}{300}$&$\frac{-1}{300}$&$\frac{3}{100}$&$\frac{-1}{150}$&$\frac{-3}{100}$&$\frac{8}{75}$&$\frac{-3}{50}$&$\frac{24}{25}$&$\frac{6}{5}$&12&$\frac{288}{5}$&1&-9&0&0&-3\\
&&&&&&&&&&&&&&&&\\
$(2233,24)$
&$\frac{1}{480}$&$\frac{-1}{480}$&$\frac{-3}{160}$&$\frac{1}{120}$&$\frac{3}{160}$&$\frac{-2}{15}$&$\frac{-3}{40}$&$\frac{6}{5}$&0&0&0&-2&-18&$\frac{3}{2}$&6&0\\
&&&&&&&&&&&&&&&&\\
$(2233,44)$ 
&$\frac{1}{1200}$&$\frac{-1}{1200}$&$\frac{3}{400}$&$\frac{-1}{150}$&$\frac{-3}{400}$&$\frac{8}{75}$&$\frac{-3}{50}$&$\frac{24}{25}$&$\frac{9}{5}$&6&$\frac{48}{5}$&$\frac{-1}{2}$&$\frac{9}{2}$&0&0&$\frac{-3}{2}$\\
&&&&&&&&&&&&&&&&\\
$(3333,11)$ 
&$\frac{1}{75}$&$\frac{-1}{100}$&$\frac{-7}{25}$&$\frac{-4}{75}$&$\frac{21}{100}$&0&$\frac{28}{25}$&0&$\frac{24}{5}$&$\frac{96}{5}$&0&0&0&4&0&0\\
&&&&&&&&&&&&&&&&\\
$(3333,12)$ 
&$\frac{1}{300}$&$\frac{-1}{200}$&$\frac{13}{100}$&$\frac{2}{75}$&$\frac{-39}{200}$&0&$\frac{26}{25}$&0&$\frac{16}{5}$&$\frac{32}{5}$&0&0&0&2&0&0\\
&&&&&&&&&&&&&&&&\\
$(3333,14)$ 
&$\frac{1}{300}$&0&$\frac{-7}{100}$&$\frac{-4}{75}$&0&0&$\frac{28}{25}$&0&$\frac{6}{5}$&$\frac{-24}{5}$&0&0&0&4&0&0\\
&&&&&&&&&&&&&&&&\\
$(3333,22)$ 
&$\frac{1}{300}$&0&$\frac{-7}{100}$&$\frac{-4}{75}$&0&0&$\frac{28}{25}$&0&$\frac{-4}{5}$&$\frac{16}{5}$&0&0&0&0&0&0\\
&&&&&&&&&&&&&&&&\\
$(3333,24)$ 
&$\frac{1}{1200}$&$\frac{-1}{400}$&$\frac{13}{400}$&$\frac{2}{75}$&$\frac{-39}{400}$&0&$\frac{26}{25}$&0&$\frac{4}{5}$&$\frac{32}{5}$&0&0&0&-1&0&0\\
&&&&&&&&&&&&&&&&\\
$(3333,44)$ 
&$\frac{1}{1200}$&$\frac{1}{400}$&$\frac{-7}{400}$&$\frac{-4}{75}$&$\frac{-21}{400}$&0&$\frac{28}{25}$&0&$\frac{-6}{5}$&$\frac{-24}{5}$&0&0&0&1&0&0\\
&&&&&&&&&&&&&&&&\\
\hline
\end{tabular}
}
%\end{landscape}
\end{center}

\newpage
\begin{center}
{\tiny 
\textbf{Table 4.} (Theorem 2.1 (ii))
\begin{tabular}{|c|c|c|c|c|c|c|c|c|c|c|c|c|c|c|}
\hline
&&&&&&&&&&&&&&\\
$a_1a_2a_3$&$ \beta_1$ &$ \beta_2$& $ \beta_3$& $ \beta_4$&$ \beta_5$ &$ \beta_6$& $ \beta_7$& $ \beta_8$&$ \beta_9$ &$ \beta_{10}$ & $ \beta_{11}$&$ \beta_{12}$& $ \beta_{13}$&$ \beta_{14}$\\
$a_4,b_1b_2$&$~~$ &$~~$& $~~$& $~~$&$~~$ &$~~$& $~~$& $~~$&$~~$ &$ ~~$ & $~~$&$~~$& $~~$&$~~$\\
\hline 
&&&&&&&&&&&&&&\\
$1112,11$&$\frac{-26}{451}$&$\frac{108}{451}$&$\frac{6656}{451}$&$\frac{27648}{451}$&$\frac{168}{451}$&$\frac{11448}{451}$&$\frac{-2496}{451}$&$\frac{-17280}{451}$&$\frac{24}{41}$&$\frac{936}{41}$&$\frac{144}{41}$&$\frac{-384}{41}$&$\frac{4032}{41}$&$\frac{-48}{41}$\\ 
&&&&&&&&&&&&&&\\
$1112,12$&$\frac{28}{451}$&$\frac{54}{451}$&$\frac{3584}{451}$&$\frac{-6912}{451}$&$\frac{480}{451}$&$\frac{-2052}{451}$&$\frac{-2688}{451}$&$\frac{1728}{451}$&$\frac{-60}{41}$&$\frac{216}{41}$&$\frac{-108}{41}$&$\frac{-2112}{41}$&$\frac{-1440}{41}$&$\frac{288}{41}$\\ 
&&&&&&&&&&&&&&\\
$1112,14$&$\frac{-26}{451}$&$\frac{108}{451}$&$\frac{1664}{451}$&$\frac{6912}{451}$&$\frac{-912}{451}$&$\frac{3672}{451}$&0&$\frac{-6912}{451}$&$\frac{-48}{41}$&$\frac{-54}{41}$&$\frac{702}{41}$&$\frac{1632}{41}$&$\frac{576}{41}$&$\frac{-228}{41}$\\ 
&&&&&&&&&&&&&&\\
$1112,22$&$\frac{-26}{451}$&$\frac{108}{451}$&$\frac{1664}{451}$&$\frac{6912}{451}$&$\frac{-912}{451}$&$\frac{3672}{451}$&0&$\frac{-6912}{451}$&$\frac{-48}{41}$&$\frac{684}{41}$&$\frac{-36}{41}$&$\frac{-2304}{41}$&$\frac{576}{41}$&$\frac{264}{41}$\\ 
&&&&&&&&&&&&&&\\
$1112,24$&$\frac{28}{451}$&$\frac{54}{451}$&$\frac{896}{451}$&$\frac{-1728}{451}$&$\frac{-66}{41}$&$\frac{-108}{451}$&$\frac{-1344}{451}$&$\frac{-864}{451}$&$\frac{-42}{41}$&$\frac{-90}{41}$&$\frac{306}{41}$&$\frac{336}{41}$&$\frac{-576}{41}$&$\frac{-36}{41}$\\ 
&&&&&&&&&&&&&&\\
$1112,44$&$\frac{-26}{451}$&$\frac{108}{451}$&$\frac{416}{451}$&$\frac{1728}{451}$&$\frac{-1182}{451}$&$\frac{1728}{451}$&$\frac{624}{451}$&$\frac{-4320}{451}$&$\frac{-66}{41}$&$\frac{252}{41}$&$\frac{288}{41}$&$\frac{-816}{41}$&$\frac{-288}{41}$&$\frac{96}{41}$\\
&&&&&&&&&&&&&&\\ 
$1222,11$&$\frac{-26}{451}$&$\frac{108}{451}$&$\frac{3328}{451}$&$\frac{13824}{451}$&$\frac{8468}{451}$&$\frac{6264}{451}$&$\frac{-15264}{451}$&$\frac{-10368}{451}$&$\frac{-352}{41}$&$\frac{-708}{41}$&$\frac{516}{41}$&$\frac{3584}{41}$&$\frac{13536}{41}$&$\frac{-660}{41}$\\ 
&&&&&&&&&&&&&&\\
$1222,12$&$\frac{28}{451}$&$\frac{54}{451}$&$\frac{1792}{451}$&$\frac{-3456}{451}$&$\frac{9598}{451}$&$\frac{-756}{451}$&$\frac{-16224}{451}$&0&$\frac{-458}{41}$&$\frac{-1956}{41}$&$\frac{168}{41}$&$\frac{2800}{41}$&$\frac{5040}{41}$&$\frac{-420}{41}$\\
&&&&&&&&&&&&&&\\
$1222,14$&$\frac{-26}{451}$&$\frac{108}{451}$&$\frac{832}{451}$&$\frac{3456}{451}$&$\frac{-6955}{451}$&$\frac{216}{41}$&$\frac{7632}{451}$&$\frac{-5184}{451}$&$\frac{473}{41}$&$\frac{1749}{41}$&$\frac{57}{41}$&-56&-144&$\frac{357}{41}$\\ 
&&&&&&&&&&&&&&\\
$1222,22$&$\frac{-26}{451}$&$\frac{108}{451}$&$\frac{832}{451}$&$\frac{3456}{451}$&$\frac{10634}{451}$&$\frac{216}{41}$&$\frac{-14016}{451}$&$\frac{-5184}{451}$&$\frac{-634}{41}$&$\frac{-2310}{41}$&$\frac{426}{41}$&112&288&$\frac{-750}{41}$\\ 
&&&&&&&&&&&&&&\\
$1222,24$&$\frac{28}{451}$&$\frac{54}{451}$&$\frac{448}{451}$&$\frac{-864}{451}$&$\frac{-3182}{451}$&$\frac{216}{451}$&$\frac{6096}{451}$&$\frac{-1296}{451}$&$\frac{166}{41}$&$\frac{474}{41}$&$\frac{6}{41}$&$\frac{-896}{41}$&$\frac{-3384}{41}$&$\frac{156}{41}$\\ 
&&&&&&&&&&&&&&\\
$1222,44$&$\frac{-26}{451}$&$\frac{108}{451}$&$\frac{208}{451}$&$\frac{864}{451}$&$\frac{2381}{451}$&$\frac{1404}{451}$&$\frac{-2880}{451}$&$\frac{-3888}{451}$&$\frac{-151}{41}$&$\frac{-681}{41}$&$\frac{219}{41}$&$\frac{1400}{41}$&$\frac{2520}{41}$&$\frac{-219}{41}$\\ 
&&&&&&&&&&&&&&\\
$1233,11$&$\frac{10}{451}$&$\frac{72}{451}$&$\frac{2560}{451}$&$\frac{-18432}{451}$&$\frac{488}{451}$&$\frac{-6192}{451}$&$\frac{-1600}{451}$&$\frac{6912}{451}$&$\frac{-112}{41}$&$\frac{-480}{41}$&$\frac{432}{41}$&$\frac{704}{41}$&$\frac{-3456}{41}$&$\frac{-24}{41}$\\
&&&&&&&&&&&&&&\\ 
$1233,12$&$\frac{-8}{451}$&$\frac{90}{451}$&$\frac{1024}{451}$&$\frac{11520}{451}$&$\frac{-1280}{451}$&$\frac{5220}{451}$&$\frac{-256}{451}$&$\frac{-8640}{451}$&$\frac{-20}{41}$&$\frac{312}{41}$&$\frac{348}{41}$&$\frac{-512}{41}$&$\frac{1440}{41}$&$\frac{24}{41}$\\
&&&&&&&&&&&&&&\\
$1233,14$&$\frac{10}{451}$&$\frac{72}{451}$&$\frac{640}{451}$&$\frac{-4608}{451}$&$\frac{-760}{451}$&$\frac{-1008}{451}$&$\frac{-640}{451}$ &0&$\frac{-64}{41}$&$\frac{-66}{41}$&$\frac{306}{41}$&$\frac{-640}{41}$&$\frac{-1152}{41}$&$\frac{96}{41}$\\ 
&&&&&&&&&&&&&&\\
$1233,22$&$\frac{10}{451}$&$\frac{72}{451}$&$\frac{640}{451}$&$\frac{-4608}{451}$&$\frac{-760}{451}$&$\frac{-1008}{451}$&$\frac{-640}{451}$&0&$\frac{-64}{41}$&$\frac{180}{41}$&$\frac{60}{41}$&$\frac{-640}{41}$&$\frac{-1152}{41}$&$\frac{96}{41}$\\ 
&&&&&&&&&&&&&&\\
$1233,24$&$\frac{-8}{451}$&$\frac{90}{451}$&$\frac{256}{451}$&$\frac{2880}{451}$&$\frac{-1238}{451}$&$\frac{180}{41}$&$\frac{128}{451}$&$\frac{-4320}{451}$&$\frac{-50}{41}$&$\frac{330}{41}$&$\frac{150}{41}$&-16&0&$\frac{72}{41}$\\ 
&&&&&&&&&&&&&&\\
$1233,44$&$\frac{10}{451}$&$\frac{72}{451}$&$\frac{160}{451}$&$\frac{-1152}{451}$&$\frac{-1072}{451}$&$\frac{288}{451}$&$\frac{-400}{451}$&$\frac{-1728}{451}$&$\frac{-52}{41}$&$\frac{222}{41}$&$\frac{90}{41}$&$\frac{-976}{41}$&$\frac{-576}{41}$&$\frac{126}{41}$\\ 
&&&&&&&&&&&&&&\\
\hline
\end{tabular}
}
%\end{landscape}
\end{center}

\begin{center}
{\tiny 
\textbf{Table 5.} (Theorem 2.1 (iii))
\begin{tabular}{|c|c|c|c|c|c|c|c|c|c|c|c|c|c|c|c|c|}
\hline
&&&&&&&&&&&&&&&&\\
$(a_1a_2a_3a_4,b_1b_2)$&$ \gamma_1$ &$ \gamma_2$& $ \gamma_3$& $ \gamma_4$&$ \gamma_5$ &$ \gamma_6$& $ \gamma_7$& $ \gamma_8$&$ \gamma_9$ &$ \gamma_{10}$ & $ \gamma_{11}$&$ \gamma_{12}$& $ \gamma_{13}$&$ \gamma_{14}$& $ \gamma_{15}$&$ \gamma_{16}$\\
\hline 
&&&&&&&&&&&&&&&&\\
$(1113,11)$&$\frac{1}{23}$&$\frac{288}{23}$&$\frac{32}{23}$&$\frac{9}{23}$&0&0&0&0&$\frac{84}{23}$&$\frac{720}{23}$&$\frac{336}{23}$&$\frac{864}{23}$&0&0&0&0\\
&&&&&&&&&&&&&&&&\\
$(1113,12)$&$\frac{1}{23}$&$\frac{144}{23}$&$\frac{-16}{23}$&$\frac{-9}{23}$&0&0&0&0&$\frac{156}{23}$&$\frac{-48}{23}$&$\frac{-168}{23}$&$\frac{-456}{23}$&0&0&0&0\\
&&&&&&&&&&&&&&&&\\
$(1113,14)$&$\frac{1}{23}$&$\frac{72}{23}$&$\frac{8}{23}$&$\frac{9}{23}$&0&0&0&0&$\frac{186}{23}$&$\frac{600}{23}$&$\frac{-228}{23}$&$\frac{372}{23}$&0&0&0&0\\
&&&&&&&&&&&&&&&&\\
$(1113,22)$&$\frac{1}{23}$&$\frac{72}{23}$&$\frac{8}{23}$&$\frac{9}{23}$&0&0&0&0&$\frac{48}{23}$&$\frac{48}{23}$&$\frac{48}{23}$&$\frac{96}{23}$&0&0&0&0\\
&&&&&&&&&&&&&&&&\\
$(1113,24)$&$\frac{1}{23}$&$\frac{36}{23}$&$\frac{-4}{23}$&$\frac{-9}{23}$&0&0&0&0&$\frac{114}{23}$&$\frac{84}{23}$&$\frac{-120}{23}$&$\frac{-156}{23}$&0&0&0&0\\
&&&&&&&&&&&&&&&&\\
$(1113,44)$&$\frac{1}{23}$&$\frac{18}{23}$&$\frac{2}{23}$&$\frac{9}{23}$&0&0&0&0&$\frac{108}{23}$&$\frac{156}{23}$&$\frac{-162}{23}$&$\frac{42}{23}$&0&0&0&0\\
&&&&&&&&&&&&&&&&\\
$(1223,11)$&0&$\frac{144}{23}$&$\frac{16}{23}$&0&$\frac{1}{23}$&0&0&$\frac{-9}{23}$&$\frac{162}{23}$&$\frac{264}{23}$&$\frac{-1188}{23}$&$\frac{420}{23}$&$\frac{192}{23}$&$\frac{864}{23}$&$\frac{-4704}{23}$&$\frac{-5760}{23}$\\
&&&&&&&&&&&&&&&&\\
$(1223,12)$&0&$\frac{72}{23}$&$\frac{-8}{23}$&0&$\frac{1}{23}$&0&0&$\frac{9}{23}$&$\frac{120}{23}$&$\frac{96}{23}$&$\frac{336}{23}$&$\frac{-384}{23}$&$\frac{-240}{23}$&$\frac{-912}{23}$&$\frac{4368}{23}$&$\frac{2784}{23}$\\
&&&&&&&&&&&&&&&&\\
$(1223,14)$&0&$\frac{36}{23}$&$\frac{4}{23}$&0&$\frac{1}{23}$&0&0&$\frac{-9}{23}$&$\frac{144}{23}$&$\frac{480}{23}$&$\frac{-504}{23}$&$\frac{312}{23}$&$\frac{-360}{23}$&$\frac{1416}{23}$&$\frac{-1944}{23}$&$\frac{-1344}{23}$\\
&&&&&&&&&&&&&&&&\\
$(1223,22)$&0&$\frac{36}{23}$&$\frac{4}{23}$&0&$\frac{1}{23}$&0&0&$\frac{-9}{23}$&$\frac{6}{23}$&$\frac{-72}{23}$&$\frac{-228}{23}$&$\frac{36}{23}$&$\frac{192}{23}$&$\frac{-240}{23}$&$\frac{-1392}{23}$&$\frac{-1344}{23}$\\
&&&&&&&&&&&&&&&&\\
$(1223,24)$&0&$\frac{18}{23}$&$\frac{-2}{23}$&0&$\frac{1}{23}$&0&0&$\frac{9}{23}$&$\frac{30}{23}$&$\frac{24}{23}$&$\frac{84}{23}$&$\frac{-96}{23}$&$\frac{36}{23}$&$\frac{-360}{23}$&$\frac{504}{23}$&$\frac{576}{23}$\\
&&&&&&&&&&&&&&&&\\
$(1223,44)$&0&$\frac{9}{23}$&$\frac{1}{23}$&0&$\frac{1}{23}$&0&0&$\frac{-9}{23}$&$\frac{36}{23}$&$\frac{120}{23}$&$\frac{-126}{23}$&$\frac{78}{23}$&$\frac{-84}{23}$&$\frac{312}{23}$&$\frac{-840}{23}$&$\frac{-240}{23}$\\
&&&&&&&&&&&&&&&&\\
$(1333,11)$&$\frac{1}{23}$&$\frac{96}{23}$&$\frac{-32}{23}$&$\frac{-3}{23}$&0&0&0&0&$\frac{260}{23}$&$\frac{352}{23}$&$\frac{-544}{23}$&$\frac{-352}{23}$&0&0&0&0\\
&&&&&&&&&&&&&&&&\\
$(1333,12)$&$\frac{1}{23}$&$\frac{48}{23}$&$\frac{16}{23}$&$\frac{3}{23}$&0&0&0&0&$\frac{116}{23}$&$\frac{160}{23}$&$\frac{-40}{23}$&$\frac{104}{23}$&0&0&0&0\\
&&&&&&&&&&&&&&&&\\
$(1333,14)$&$\frac{1}{23}$&$\frac{24}{23}$&$\frac{-8}{23}$&$\frac{-3}{23}$&0&0&0&0&$\frac{170}{23}$&$\frac{16}{23}$&$\frac{-292}{23}$&$\frac{-340}{23}$&0&0&0&0\\
&&&&&&&&&&&&&&&&\\
$(1333,22)$&$\frac{1}{23}$&$\frac{24}{23}$&$\frac{-8}{23}$&$\frac{-3}{23}$&0&0&0&0&$\frac{32}{23}$&$\frac{16}{23}$&$\frac{-16}{23}$&$\frac{-64}{23}$&0&0&0&0\\
&&&&&&&&&&&&&&&&\\
$(1333,24)$&$\frac{1}{23}$&$\frac{12}{23}$&$\frac{4}{23}$&$\frac{3}{23}$&0&0&0&0&$\frac{26}{23}$&$\frac{76}{23}$&$\frac{32}{23}$&$\frac{116}{23}$&0&0&0&0\\
&&&&&&&&&&&&&&&&\\
$(1333,44)$&$\frac{1}{23}$&$\frac{6}{23}$&$\frac{-2}{23}$&$\frac{-3}{23}$&0&0&0&0&$\frac{44}{23}$&$\frac{-68}{23}$&$\frac{-22}{23}$&$\frac{-130}{23}$&0&0&0&0\\
&&&&&&&&&&&&&&&&\\
\hline
\end{tabular}
}
%\end{landscape}
\end{center}

\begin{center}
{\tiny 
\textbf{Table 6.} (Theorem 2.1 (iv))
\begin{tabular}{|c|c|c|c|c|c|c|c|c|c|c|c|c|c|c|}
\hline
&&&&&&&&&&&&&&\\
$a_1a_2a_3$&$ \delta_1$ &$ \delta_2$& $ \delta_3$& $ \delta_4$&$ \delta_5$ &$ \delta_6$& $ \delta_7$& $ \delta_8$&$ \delta_9$ &$ \delta_{10}$ & $ \delta_{11}$&$ \delta_{12}$& $ \delta_{13}$&$ \delta_{14}$\\
$a_4b_1b_2$&$~~$ &$~~$& $~~$& $~~$&$~~$ &$~~$& $~~$& $~~$&$~~$ &$ ~~$ & $~~$&$~~$& $~~$&$~~$\\
\hline 
&&&&&&&&&&&&&&\\
$112311$&$\frac{1}{261}$&$\frac{256}{29}$&$\frac{-256}{261}$&$\frac{-1}{29}$&$\frac{1808}{87}$&$\frac{656}{29}$&$\frac{-2056}{87}$&$\frac{-3808}{29}$&$\frac{-4144}{29}$&$\frac{736}{3}$&$\frac{472}{3}$&$\frac{-41984}{87}$&$\frac{-1096}{87}$&$\frac{-968}{87}$\\
&&&&&&&&&&&&&&\\
$112312$&$\frac{1}{261}$&$\frac{128}{29}$&$\frac{128}{261}$&$\frac{1}{29}$&$\frac{208}{87}$&$\frac{-32}{29}$&$\frac{-284}{87}$&$\frac{-368}{29}$&$\frac{1048}{29}$&$\frac{-6224}{87}$&$\frac{-7100}{87}$&$\frac{21248}{87}$&$\frac{8}{3}$&$\frac{500}{87}$\\
&&&&&&&&&&&&&&\\ 
$112314$&$\frac{1}{261}$&$\frac{64}{29}$&$\frac{-64}{261}$&$\frac{-1}{29}$&$\frac{-84}{29}$&$\frac{-114}{29}$&$\frac{450}{29}$&$\frac{1800}{29}$&$\frac{-324}{29}$&$\frac{-4200}{29}$&$\frac{-2814}{29}$&$\frac{-3072}{29}$&$\frac{318}{29}$&$\frac{414}{29}$\\
&&&&&&&&&&&&&&\\ 
$112322$&$\frac{1}{261}$&$\frac{64}{29}$&$\frac{-64}{261}$&$\frac{-1}{29}$&$\frac{264}{29}$&$\frac{60}{29}$&$\frac{-420}{29}$&$\frac{-1680}{29}$&$\frac{-1368}{29}$&$\frac{2064}{29}$&$\frac{1884}{29}$&$\frac{-3072}{29}$&$\frac{-204}{29}$&$\frac{-108}{29}$\\
&&&&&&&&&&&&&&\\
$112324$&$\frac{1}{261}$&$\frac{32}{29}$&$\frac{32}{261}$&$\frac{1}{29}$&$\frac{860}{87}$&$\frac{218}{29}$&$\frac{-970}{87}$&$\frac{-2296}{29}$&$\frac{-628}{29}$&$\frac{7976}{87}$&$\frac{-778}{87}$&$\frac{4288}{87}$&$\frac{-622}{87}$&$\frac{-10}{3}$\\
&&&&&&&&&&&&&&\\ 
$112344$&$\frac{1}{261}$&$\frac{16}{29}$&$\frac{-16}{261}$&$\frac{-1}{29}$&$\frac{16}{87}$&$\frac{-176}{29}$&$\frac{244}{87}$&$\frac{592}{29}$&$\frac{-152}{29}$&$\frac{-6992}{87}$&$\frac{-3404}{87}$&$\frac{-1024}{87}$&$\frac{292}{87}$&$\frac{620}{87}$\\
&&&&&&&&&&&&&&\\ 
$222311$&$\frac{1}{261}$&$\frac{128}{29}$&$\frac{-128}{261}$&$\frac{-1}{29}$&$\frac{5480}{261}$&$\frac{1472}{87}$&$\frac{-15016}{261}$&$\frac{-2992}{87}$&$\frac{-23584}{87}$&$\frac{79664}{261}$&$\frac{102248}{261}$&$\frac{-194048}{261}$&$\frac{-116}{9}$&$\frac{-3704}{261}$\\
&&&&&&&&&&&&&&\\ 
$222312$&$\frac{1}{261}$&$\frac{64}{29}$&$\frac{64}{261}$&$\frac{1}{29}$&$\frac{-160}{261}$&$\frac{-640}{87}$&$\frac{7172}{261}$&$\frac{-5656}{87}$&$\frac{12320}{87}$&$\frac{-50824}{261}$&$\frac{-47284}{261}$&$\frac{96640}{261}$&$\frac{1076}{261}$&$\frac{964}{261}$\\
&&&&&&&&&&&&&&\\ 
$222314$&$\frac{1}{261}$&$\frac{32}{29}$&$\frac{-32}{261}$&$\frac{-1}{29}$&$\frac{824}{261}$&$\frac{-466}{87}$&$\frac{4970}{261}$&$\frac{-1888}{87}$&$\frac{-628}{87}$&$\frac{-24496}{261}$&$\frac{-10030}{261}$&$\frac{-44672}{261}$&$\frac{494}{261}$&$\frac{394}{261}$\\
&&&&&&&&&&&&&&\\ 
$222322$&$\frac{1}{261}$&$\frac{32}{29}$&$\frac{-32}{261}$&$\frac{-1}{29}$&$\frac{824}{261}$&$\frac{1100}{87}$&$\frac{-9124}{261}$&$\frac{2288}{87}$&$\frac{-10024}{87}$&$\frac{50672}{261}$&$\frac{32252}{261}$&$\frac{-44672}{261}$&$\frac{-1072}{261}$&$\frac{-1172}{261}$\\
&&&&&&&&&&&&&&\\ 
$222324$&$\frac{1}{261}$&$\frac{16}{29}$&$\frac{16}{261}$&$\frac{1}{29}$&$\frac{-748}{261}$&$\frac{518}{87}$&$\frac{-2470}{261}$&$\frac{2936}{87}$&$\frac{-1156}{87}$&$\frac{11192}{261}$&$\frac{-8830}{261}$&$\frac{21088}{261}$&$\frac{578}{261}$&$\frac{562}{261}$\\
&&&&&&&&&&&&&&\\ 
$222344$&$\frac{1}{261}$&$\frac{8}{29}$&$\frac{-8}{261}$&$\frac{-1}{29}$&$\frac{-340}{261}$&$\frac{224}{87}$&$\frac{-604}{261}$&$\frac{1520}{87}$&$\frac{-1936}{87}$&$\frac{5840}{261}$&$\frac{-6388}{261}$&$\frac{-7328}{261}$&$\frac{284}{261}$&$\frac{244}{261}$\\
&&&&&&&&&&&&&&\\ 
$233311$&$\frac{1}{261}$&$\frac{256}{87}$&$\frac{256}{261}$&$\frac{1}{87}$&$\frac{-13360}{261}$&$\frac{-3520}{87}$&$\frac{20816}{261}$&$\frac{12608}{29}$&$\frac{20960}{87}$&$\frac{-219712}{261}$&$\frac{-113968}{261}$&$\frac{133120}{261}$&$\frac{15464}{261}$&$\frac{16384}{261}$\\
&&&&&&&&&&&&&&\\ 
$233312$&$\frac{1}{261}$&$\frac{128}{87}$&$\frac{-128}{261}$&$\frac{-1}{87}$&$\frac{11168}{261}$&$\frac{2312}{87}$&$\frac{-14212}{261}$&$\frac{-27920}{87}$&$\frac{-14968}{87}$&$\frac{129296}{261}$&$\frac{45212}{261}$&$\frac{-64256}{261}$&$\frac{-9856}{261}$&$\frac{-8948}{261}$\\
&&&&&&&&&&&&&&\\ 
$233314$&$\frac{1}{261}$&$\frac{64}{87}$&$\frac{64}{261}$&$\frac{1}{87}$&$\frac{3836}{261}$&$\frac{338}{87}$&$\frac{-4126}{261}$&$\frac{-2968}{29}$&$\frac{-2092}{87}$&$\frac{21368}{261}$&$\frac{-10174}{261}$&$\frac{32512}{261}$&$\frac{-2530}{261}$&$\frac{-1514}{261}$\\
&&&&&&&&&&&&&&\\ 
$233322$&$\frac{1}{261}$&$\frac{64}{87}$&$\frac{64}{261}$&$\frac{1}{87}$&$\frac{-5560}{261}$&$\frac{-1228}{87}$&$\frac{6836}{261}$&$\frac{4688}{29}$&$\frac{7304}{87}$&$\frac{-72592}{261}$&$\frac{-24268}{261}$&$\frac{32512}{261}$&$\frac{5300}{261}$&$\frac{6316}{261}$\\
&&&&&&&&&&&&&&\\ 
$233324$&$\frac{1}{261}$&$\frac{32}{87}$&$\frac{-32}{261}$&$\frac{-1}{87}$&$\frac{-940}{261}$&$\frac{-274}{87}$&$\frac{842}{261}$&$\frac{2584}{87}$&$\frac{-76}{87}$&$\frac{-12808}{261}$&$\frac{2666}{261}$&$\frac{-14528}{261}$&$\frac{878}{261}$&$\frac{1882}{261}$\\ 
&&&&&&&&&&&&&&\\
$233344$&$\frac{1}{261}$&$\frac{16}{87}$&$\frac{16}{261}$&$\frac{1}{87}$&$\frac{1088}{261}$&$\frac{128}{87}$&$\frac{-2140}{261}$&$\frac{-1120}{29}$&$\frac{-808}{87}$&$\frac{11168}{261}$&$\frac{5204}{261}$&$\frac{7360}{261}$&$\frac{-1156}{261}$&$\frac{-4}{9}$\\
&&&&&&&&&&&&&&\\
\hline
\end{tabular}
}

\end{center}
%\end{landscape}

\begin{center}
{\tiny 
\textbf{Table 7.} (Theorem 2.2)
\begin{tabular}{|c|c|c|c|c|c|c|c|c|c|c|c|c|c|c|c|c|}
\hline
&&&&&&&&&&&&&&&&\\
$(1,c_1,c_2,c_3)$&$ \nu_1$ &$ \nu_2$& $ \nu_3$& $ \nu_4$&$ \nu_5$ &$ \nu_6$& $ \nu_7$& $ \nu_8$&$ \nu_9$ &$ \nu_{10}$ & $ \nu_{11}$&$ \nu_{12}$& $ \nu_{13}$&$ \nu_{14}$& $ \nu_{15}$&$ \nu_{16}$\\
\hline 
&&&&&&&&&&&&&&&&\\
$(1,1,1,2)$&$\frac{3}{40}$&$\frac{-1}{5}$&$\frac{-27}{40}$&0&$\frac{9}{5}$&0&0&0&0&0&0&0&0&0&0&0\\
&&&&&&&&&&&&&&&&\\
$(1,1,1,4)$&$\frac{3}{100}$&$\frac{-9}{100}$&$\frac{27}{100}$&$\frac{4}{25}$&$\frac{-81}{100}$&0&$\frac{36}{25}$&0&$\frac{54}{5}$&$\frac{432}{5}$&0&0&0&0&0&0\\
&&&&&&&&&&&&&&&&\\
$(1,1,1,8)$&$\frac{3}{160}$&$\frac{-9}{160}$&$\frac{-27}{160}$&$\frac{9}{80}$&$\frac{81}{160}$&$\frac{-1}{5}$&$\frac{-81}{80}$&$\frac{9}{5}$&0&0&0&$\frac{-27}{4}$&$\frac{243}{4}$&0&81&$\frac{81}{4}$\\
&&&&&&&&&&&&&&&&\\
$(1,1,2,2)$&$\frac{1}{50}$&$\frac{2}{25}$&$\frac{9}{50}$&0&$\frac{18}{25}$&0&0&0&$\frac{36}{5}$&0&0&0&0&0&0&0\\
&&&&&&&&&&&&&&&&\\
$(1,1,2,4)$&$\frac{1}{80}$&$\frac{1}{16}$&$\frac{-9}{80}$&$\frac{-1}{5}$&$\frac{-9}{16}$&0&$\frac{9}{5}$&0&0&0&0&0&0&9&0&0\\
&&&&&&&&&&&&&&&&\\
$(1,1,2,8)$&$\frac{1}{200}$&$\frac{1}{40}$&$\frac{9}{200}$&$\frac{-9}{100}$&$\frac{9}{40}$&$\frac{4}{25}$&$\frac{-81}{100}$&$\frac{36}{25}$&$\frac{9}{5}$&$\frac{144}{5}$&$\frac{1152}{5}$&9&81&0&0&0\\
&&&&&&&&&&&&&&&&\\
$(1,1,4,4)$&$\frac{1}{200}$&$\frac{3}{200}$&$\frac{9}{200}$&$\frac{2}{25}$&$\frac{27}{200}$&0&$\frac{18}{25}$&0&$\frac{54}{5}$&$\frac{216}{5}$&0&0&0&0&0&0\\
&&&&&&&&&&&&&&&&\\
$(1,1,4,8)$&$\frac{1}{320}$&$\frac{3}{320}$&$\frac{-9}{320}$&$\frac{1}{16}$&$\frac{-27}{320}$&$\frac{-1}{5}$&$\frac{-9}{16}$&$\frac{9}{5}$&0&0&0&$\frac{9}{4}$&$\frac{-81}{4}$&$\frac{9}{4}$&27&$\frac{27}{4}$\\
&&&&&&&&&&&&&&&&\\
$(1,1,8,8)$&$\frac{1}{800}$&$\frac{3}{800}$&$\frac{9}{800}$&$\frac{3}{200}$&$\frac{27}{800}$&$\frac{2}{25}$&$\frac{27}{200}$&$\frac{18}{25}$&$\frac{36}{5}$&$\frac{234}{5}$&$\frac{576}{5}$&$\frac{9}{2}$&$\frac{81}{2}$&0&0&0\\
&&&&&&&&&&&&&&&&\\
$(1,2,2,2)$&$\frac{1}{40}$&$\frac{-3}{20}$&$\frac{-9}{40}$&0&$\frac{27}{20}$&0&0&0&0&0&0&0&0&0&0&0\\
&&&&&&&&&&&&&&&&\\
$(1,2,2,4)$&$\frac{1}{100}$&$\frac{-7}{100}$&$\frac{9}{100}$&$\frac{4}{25}$&$\frac{-63}{100}$&0&$\frac{36}{25}$&0&$\frac{18}{5}$&$\frac{72}{5}$&0&0&0&0&0&0\\
&&&&&&&&&&&&&&&&\\
$(1,2,2,8)$&$\frac{1}{160}$&$\frac{-7}{160}$&$\frac{-9}{160}$&$\frac{9}{80}$&$\frac{63}{160}$&$\frac{-1}{5}$&$\frac{-81}{80}$&$\frac{9}{5}$&0&0&0&$\frac{-9}{4}$&$\frac{81}{4}$&0&9&$\frac{27}{4}$\\
&&&&&&&&&&&&&&&&\\
$(1,2,4,4)$&$\frac{1}{160}$&$\frac{-1}{32}$&$\frac{-9}{160}$&$\frac{-1}{10}$&$\frac{9}{32}$&0&$\frac{9}{10}$&0&0&0&0&0&0&$\frac{9}{2}$&0&0\\
&&&&&&&&&&&&&&&&\\
$(1,2,4,8)$&$\frac{1}{400}$&$\frac{-1}{80}$&$\frac{9}{400}$&$\frac{-1}{20}$&$\frac{-9}{80}$&$\frac{4}{25}$&$\frac{-9}{20}$&$\frac{36}{25}$&$\frac{9}{10}$&$\frac{27}{5}$&$\frac{72}{5}$&$\frac{9}{2}$&$\frac{81}{2}$&0&0&0\\
&&&&&&&&&&&&&&&&\\
$(1,2,8,8)$&$\frac{1}{640}$&$\frac{-1}{128}$&$\frac{-9}{640}$&$\frac{-3}{160}$&$\frac{9}{128}$&$\frac{-1}{10}$&$\frac{27}{160}$&$\frac{9}{10}$&0&0&0&$\frac{-9}{8}$&$\frac{81}{8}$&$\frac{27}{8}$&$\frac{9}{2}$&$\frac{27}{8}$\\
&&&&&&&&&&&&&&&&\\
$(1,4,4,4)$&$\frac{1}{400}$&$\frac{-9}{400}$&$\frac{9}{400}$&$\frac{3}{25}$&$\frac{-81}{400}$&0&$\frac{27}{25}$&0&$\frac{27}{5}$&$\frac{54}{5}$&0&0&0&0&0&0\\
&&&&&&&&&&&&&&&&\\
$(1,4,4,8)$&$\frac{1}{640}$&$\frac{-9}{640}$&$\frac{-9}{640}$&$\frac{7}{80}$&$\frac{81}{640}$&$\frac{-1}{5}$&$\frac{-63}{80}$&$\frac{9}{5}$&0&0&0&$\frac{9}{8}$&$\frac{-81}{8}$&$\frac{9}{8}$&0&$\frac{27}{8}$\\
&&&&&&&&&&&&&&&&\\
$(1,4,8,8)$&$\frac{1}{1600}$&$\frac{-9}{1600}$&$\frac{9}{1600}$&$\frac{1}{40}$&$\frac{-81}{1600}$&$\frac{2}{25}$&$\frac{9}{40}$&$\frac{18}{25}$&$\frac{18}{5}$&$\frac{36}{5}$&$\frac{36}{5}$&$\frac{9}{4}$&$\frac{81}{4}$&0&0&0\\
&&&&&&&&&&&&&&&&\\
$(1,8,8,8)$&$\frac{1}{2560}$&$\frac{-9}{2560}$&$\frac{-9}{2560}$&$\frac{9}{320}$&$\frac{81}{2560}$&$\frac{-3}{20}$&$\frac{-81}{320}$&$\frac{27}{20}$&0&0&0&$\frac{27}{32}$&$\frac{-243}{32}$&$\frac{81}{32}$&0&$\frac{81}{32}$\\
&&&&&&&&&&&&&&&&\\
\hline
\end{tabular}
}
%\end{landscape}
\end{center}

%\smallskip

\vglue -0.7cm

\subsection{Sample formulas}
In this section we shall give explicit formulas for a few cases from Tables 1 and 2. 

%For Table 1, we give formulas for the cases (1,1,1,1,1,1), (1,1,1,1,1,2)(\thmref{1}(i)); %(1,1,1,2,1,1),(1,1,1,2,1,2) (\thmref{1}(ii));(1,1,1,3,1,1),\\(1,1,1,3,1,2) (\thmref{1}(iii)); %(1,1,2,3,1,1),(1,1,2,3,1,2),(\thmref{1}(iv)), and  for Table 2, the formulas are given for the cases %(1,1,1,2), (1,1,1,4) (\thmref{2}). 

%\smallskip
 
\noindent {\bf First two formulas of \thmref{1}(i):}
\begin{equation*}
\begin{split}
%{\rm (i)} ~
{\bf N(1,1,1,1,1,1;n)}&=\frac{112}{5}\sigma_{3}(n)-\frac{84}{5} \sigma_{3}(n/2)- \frac{432}{5} \sigma_{3}(n/3)-\frac{448}{5} \sigma_{3}(n/4) +\frac{324}{5}\sigma_{3}(n/6)\\
&+ \frac{1728}{5} \sigma_{3}(n/12)-\frac{72}{5} a_{4,6}(n)-\frac{288}{5}a_{4,6}(n/2)+12 
a_{4,12}(n),\\
\end{split}
\end{equation*}
%\smallskip
\begin{equation*}
\begin{split}
%{\rm (i)} ~
{\bf N(1,1,1,1,1,2;n)}&=\frac{52}{5}\sigma_{3}(n)-\frac{78}{5} \sigma_{3}(n/2)+ \frac{108}{5} \sigma_{3}(n/3)+\frac{416}{5} \sigma_{3}(n/4) -\frac{324}{5}\sigma_{3}(n/6)\\
&+ \frac{864}{5} \sigma_{3}(n/12)+\frac{48}{5} a_{4,6}(n)+\frac{96}{5}a_{4,6}(n/2)-6 
a_{4,12}(n).\\
\end{split}
\end{equation*}
%\end{thm}
%}
%\smallskip
\noindent {\bf First two formulas of \thmref{1}(ii):}
\begin{equation*}
\begin{split}
%{\rm (i)} ~
{\bf N(1,1,1,2,1,1;n)}&=-\frac{26}{451} \sigma_{3;{\bf 1},\chi_{8}}(n)+\frac{108}{451} \sigma_{3;{\bf 1},\chi_{8}}(n/3)+\frac{6656}{451} \sigma_{3;\chi_{8},{\bf 1}}(n)+ \frac{27648}{451} \sigma_{3;\chi_{8},{\bf 1}}(n/3)\\
& + \frac{168}{451} a_{4,8,\chi_8;1}(n)+\frac{11448}{451} a_{4,8,\chi_8;1}(n/3)- \frac{2496}{451} a_{4,8,\chi_8;2}(n)-\frac{17280}{451} a_{4,8,\chi_8;2}(n/3)\\
&+ \frac{24}{41} a_{4,24,\chi_8;1}(n)+\frac{936}{41} a_{4,24,\chi_8;2}(n)+\frac{144}{41} a_{4,24,\chi_8;3}(n)- \frac{384}{41} a_{4,24,\chi_8;4}(n)\\
&+ \frac{4032}{41} a_{4,24,\chi_8;5}(n)-\frac{48}{41} a_{4,24,\chi_8;6}(n),\\
\end{split}
\end{equation*}

\begin{equation*}
\begin{split}
%{\rm (i)} ~
{\bf N(1,1,1,2,1,2;n)}&=\frac{28}{451} \sigma_{3;{\bf 1},\chi_{8}}(n)+\frac{54}{451} \sigma_{3;{\bf 1},\chi_{8}}(n/3)+\frac{3584}{451} \sigma_{3;\chi_{8},{\bf 1}}(n)- \frac{6912}{451} \sigma_{3;\chi_{8},{\bf 1}}(n/3)\\
& +\frac{480}{451} a_{4,8,\chi_8;1}(n)-\frac{2052}{451} a_{4,8,\chi_8;1}(n/3)-\frac{2688}{451} a_{4,8,\chi_8;2}(n)+\frac{1728}{451} a_{4,8,\chi_8;2}(n/3)\\
&- \frac{60}{41} a_{4,24,\chi_8;1}(n)+\frac{216}{41} a_{4,24,\chi_8;2}(n)-\frac{108}{41} a_{4,24,\chi_8;3}(n)- \frac{2112}{41} a_{4,24,\chi_8;4}(n)\\
&-\frac{1440}{41} a_{4,24,\chi_8;5}(n)+\frac{288}{41} a_{4,24,\chi_8;6}(n).\\
\end{split}
\end{equation*}
%For the space \textbf{\textbf{$M_4({\Gamma_0(48)})$}} (1,2,3,4,6) case

\smallskip

\noindent {\bf First two formulas of \thmref{1}(iii):}
\begin{equation*}
\begin{split}
%{\rm (i)} ~
{\bf N(1,1,1,3,1,1;n)}&=\frac{1}{23} \sigma_{3;{\bf 1},\chi_{12},}(n)+\frac{288}{23} \sigma_{3;\chi_{12},{\bf 1},}(n)+\frac{32}{23} \sigma_{3;\chi_{-4},\chi_{-3}}(n)+\frac{9}{23} \sigma_{3;\chi_{-3},\chi_{-4}}(n)\\
&+ \frac{84}{23}a_{4,12,\chi_{12};1}(n)+\frac{720}{23} a_{4,12,\chi_{12};2}(n)+ \frac{336}{23}a_{4,12,\chi_{12};3}(n)+\frac{864}{23}a_{4,12,\chi_{12};4}(n),\\
\end{split}
\end{equation*}
\begin{equation*}
\begin{split}
%{\rm (i)} ~
{\bf N(1,1,1,3,1,2;n)}&=\frac{1}{23} \sigma_{3;{\bf 1},\chi_{12},}(n)+\frac{144}{23} \sigma_{3;\chi_{12},{\bf 1},}(n)-\frac{16}{23} \sigma_{3;\chi_{-4},\chi_{-3}}(n)-\frac{9}{23} \sigma_{3;\chi_{-3},\chi_{-4}}(n)\\
&+ \frac{156}{23}a_{4,12,\chi_{12};1}(n)-\frac{48}{23} a_{4,12,\chi_{12};2}(n)- \frac{168}{23}a_{4,12,\chi_{12};3}(n)-\frac{456}{23}a_{4,12,\chi_{12};4}(n).\\
\end{split}
\end{equation*}
%{\tiny

\smallskip

\noindent {\bf First two formulas of \thmref{1}(iv):}
\begin{equation*}
\begin{split}
%{\rm (i)} ~
{\bf N(1,1,2,3,1,1;n)}&=\frac{1}{261} \sigma_{3;{\bf 1},\chi_{24}}(n)+ \frac{256}{29} \sigma_{3;\chi_{24},{\bf 1}}(n)-\frac{256}{261} \sigma_{3;\chi_{-8},\chi_{-3}}(n)-\frac{1}{29} \sigma_{3;\chi_{-3},\chi_{-8}}(n)\\
&+ \frac{1808}{87} a_{4,24,\chi_{24};1}(n)+\frac{656}{29} a_{4,24,\chi_{24};2}(n)-\frac{2056}{87} a_{4,24,\chi_{24};3}(n)-\frac{3808}{29} a_{4,24,\chi_{24};4}(n)\\
& - \frac{4144}{29} a_{4,24,\chi_{24};5}(n)+\frac{736}{3} a_{4,24,\chi_{24};6}(n)+ \frac{472}{3} a_{4,24,\chi_{24};7}(n)-\frac{41984}{87} a_{4,24,\chi_{24};8}(n)\\
&- \frac{1096}{87} a_{4,24,\chi_{24};9}(n)-\frac{968}{87} a_{4,24,\chi_{24};10}(n),\\
\end{split}
\end{equation*}

\begin{equation*}
\begin{split}
%{\rm (i)} ~
{\bf N(1,1,2,3,1,2;n)}&=\frac{1}{261} \sigma_{3;{\bf 1},\chi_{24}}(n)+ \frac{128}{29} \sigma_{3;\chi_{24},{\bf 1}}(n)+\frac{128}{261} \sigma_{3;\chi_{-8},\chi_{-3}}(n)+\frac{1}{29} \sigma_{3;\chi_{-3},\chi_{-8}}(n)\\
&+\frac{208}{87} a_{4,24,\chi_{24};1}(n)-\frac{32}{29} a_{4,24,\chi_{24};2}(n)-\frac{284}{87} a_{4,24,\chi_{24};3}(n)-\frac{368}{29} a_{4,24,\chi_{24};4}(n)\\
&+ \frac{1048}{29} a_{4,24,\chi_{24};5}(n)-\frac{6224}{87} a_{4,24,\chi_{24};6}(n)-\frac{7100}{87} a_{4,24,\chi_{24};7}(n)+\frac{21248}{87} a_{4,24,\chi_{24};8}(n)\\
&+\frac{8}{3} a_{4,24,\chi_{24};9}(n)+\frac{500}{87} a_{4,24,\chi_{24};10}(n).\\
\end{split}
\end{equation*}
%\begin{thm}\label{1}
%Let $n\in {\mathbb N}$. Then 
%\noindent  {\bf The cases (0,0,0,2,6) and (0,0,0,6,2)  in (Table 2)}

\smallskip

\noindent {\bf First two formulas of \thmref{2}:}
\begin{equation*}
\begin{split}
%{\rm (i)} ~
{\bf M(1,1,1,2;n)}& = 18 \sigma_{3}(n)- 48 \sigma_{3}(n/2)-  162 \sigma_{3}(n/3) + 432 \sigma_{3}(n/6),\\
%\end{split}
%\end{equation*}
%\begin{equation*}
%\begin{split}
%{\rm (i)} ~
{\bf M(1,1,1,4;n)}&=\frac{36}{5}\sigma_{3}(n)- 48 \sigma_{3}(n/2)+\frac{324}{5} \sigma_{3}(n/3)+\frac{192}{5} \sigma_{3}(n/4) -\frac{972}{5}\sigma_{3}(n/6)\\
&+\frac{1728}{5} \sigma_{3}(n/12)+\frac{54}{5} a_{4,6}(n)+\frac{432}{5}a_{4,6}(n/2).\\
\end{split}
\end{equation*}

\smallskip

\bigskip

{\small 
\noindent {\bf Acknowledgements.} ~
We have used the open-source mathematics software SAGE (www.sagemath.org) to perform our calculations. Brundaban Sahu is partially funded by SERB grant SR/FTP/MS-053/2012. 
Anup Kumar Singh thanks the Department of Mathematics, NISER, Bhubaneswar (where part of this work is done) for their warm hospitality and support. 
}

\end{document}